\documentclass[10pt]{article}
\usepackage{geometry}
\usepackage{floatrow}
\usepackage{graphicx}
\usepackage{amsmath,amssymb}
\usepackage{mathtools}
\usepackage{bm}
\usepackage{booktabs}
\usepackage{array}
\usepackage{subcaption}
\usepackage{siunitx}
\usepackage{hyperref}
\usepackage{algorithm}
\usepackage[noend]{algpseudocode}
\usepackage{listings,xcolor}
\usepackage{varwidth}
\usepackage{authblk}

\graphicspath{{figures/}}

\newfloatcommand{capbtabbox}{table}[][\FBwidth]
\let\Vec\undefined
\DeclareMathOperator{\Vec}{vec}
\DeclareMathOperator\erf{erf}

\algrenewcommand{\algorithmiccomment}[1]
{\hfill\hphantom{}\hspace{0.01em}\hspace*{\fill}\mbox{$\rhd$ \textit{#1}}}

\title{Approximate tensor-product preconditioners for very high order 
       discontinuous Galerkin methods}
\author[1]{Will Pazner}
\author[2]{Per-Olof Persson}
\affil[1]{Division of Applied Mathematics, 
          Brown University, Providence, RI, 02912}
\affil[2]{Department of Mathematics, 
          University of California, Berkeley, Berkeley, CA, 94720-3840}

\date{}

\newcommand{\keywords}[1]{ {\small \textbf{\textit{Keywords:}} #1} }

\begin{document}
\maketitle
\begin{abstract}
In this paper, we develop a new tensor-product based preconditioner for
discontinuous Galerkin methods with polynomial degrees higher than those
typically employed. This preconditioner uses an automatic, purely algebraic
method to approximate the exact block Jacobi preconditioner by Kronecker
products of several small, one-dimensional matrices. Traditional matrix-based
preconditioners require $\mathcal{O}(p^{2d})$ storage and $\mathcal{O}(p^{3d})$
computational work, where $p$ is the degree of basis polynomials used, and $d$
is the spatial dimension. Our SVD-based tensor-product preconditioner requires
$\mathcal{O}(p^{d+1})$ storage, $\mathcal{O}(p^{d+1})$ work in two spatial
dimensions, and $\mathcal{O}(p^{d+2})$ work in three spatial dimensions.
Combined with a matrix-free Newton-Krylov solver, these preconditioners allow
for the solution of DG systems in linear time in $p$ per degree of freedom in
2D, and reduce the computational complexity from $\mathcal{O}(p^9)$ to
$\mathcal{O}(p^5)$ in 3D. Numerical results are shown in 2D and 3D for the
advection, Euler, and Navier-Stokes equations, using polynomials of degree up to
$p=30$. For many test cases, the preconditioner results in similar iteration 
counts when compared with the exact block Jacobi preconditioner, and performance
is significantly improved for high polynomial degrees $p$.
\end{abstract}

%%%%%%%%%%%%%%%%%%%%%%%%%%%%%%%%%%%%%%%%%%%%%%%%%%%%%%%%%%%%%%%%%%%%%%%%%%%%%%%%
% Introduction
%%%%%%%%%%%%%%%%%%%%%%%%%%%%%%%%%%%%%%%%%%%%%%%%%%%%%%%%%%%%%%%%%%%%%%%%%%%%%%%%

\section{Introduction}\label{sec:intro}

The discontinuous Galerkin (DG) method, introduced in \cite{Reed:1973} by Reed
and Hill for the neutron transport equation, is a finite element method using
discontinuous basis functions. In the 1990s, the DG method was extended to
nonlinear systems of conservation laws by Cockburn and Shu
\cite{Cockburn:1998:RKDG}. The method has many attractive features, including
arbitrarily high formal order of accuracy, and the ability to use general,
unstructured meshes with complex geometry. In particular, the promise of a
high-order method for fluid flow problems has spurred recent recent interest in
the DG method \cite{Peraire:2011}. Higher-order methods promise highly-accurate
solutions for less computational cost than traditional low-order methods.
Additionally, high-order methods are more computationally intensive per degree
of freedom than corresponding low-order methods, resulting in a higher
computation-to-communication ratio, and thus rendering these method more
amenable to parallelization \cite{baggag1999parallel}.

High-order accuracy is achieved with the DG method by using a high-degree local
polynomial basis on each element in the mesh. There are several challenges that
can prevent the use of very high-degree polynomials as basis functions. The
number of degrees of freedom per element scales as $\mathcal{O}(p^d)$, where $p$
is the degree of polynomial approximation, and $d$ is the spatial dimension,
resulting in very computationally expensive methods. Using tensor-product
evaluations and sum factorizations \cite{Orszag:1980}, it is possible to reduce
the computational cost of these methods,
{%
however, the spectrum of the semi-discrete operator grows at a rate bounded
above by $(p+1)(p+2)/h$, and well approximated by $(p+1)^{1.78}/h$ where $p$ is
the degree of polynomial approximation, and $h$ is the element size
\cite{gottlieb1991cfl,Warburton2008}. As a result, when using explicit time
integration schemes, the time step must satisfy a restrictive stability
condition given by (approximately) $\Delta t \leq C h / (p+1)^{1.78}$
\cite{Krivodonova:2013}.}
On the other hand, the DG method couples all the degrees of freedom within each
element, so that implicit time integration methods result in block-structured
systems of equations, with blocks of size $p^d \times p^d$. Strategies for
solving these large linear systems include Newton-Krylov iterative solvers
coupled with an appropriate preconditioner \cite{Persson:2008}. Examples of
preconditioners considered include block Jacobi and Gauss-Seidel
\cite{Mardal:2007}, incomplete LU factorizations (LU)
\cite{Persson:2006efficient}, and domain decomposition techniques
\cite{Diosady:2011domain}. Multigrid and multi-level solvers have also been
considered \cite{Gopalakrishnan:2003, Kanschat:2008, Birken:2013}.

Many of the above preconditioners require the inversion of large the $p^d \times
p^d$ blocks corresponding to each element. Using dense linear algebra, this
requires $\mathcal{O}(p^{3d})$operations, which quickly becomes intractable. One
approach to reduce the computational complexity of implicit methods is to
combine Kronecker and sum-factorization techniques with a matrix-free approach.
Matrix-free approaches for the DG method have been considered in \textit{e.g.}
\cite{Crivellini2011} and \cite{klofkorn2015efficient}. Past work on efficiently
preconditioning these systems includes the use of alternating-direction-implicit
(ADI) and fast diagonalization method (FDM) preconditioners \cite{Diosady:2017}.
Kronecker-product approaches have been studied in the context of spectral
methods \cite{Shen2011}, and applications to the Navier-Stokes equations were
considered in \cite{EscobarVargas}. In this work, we describe a new approximate
Kronecker-product preconditioner that, when combined with a matrix-free tensor
product evaluation approach, allows for efficient solution of the linear systems
that arise from implicit time discretizations for high polynomial degree DG
methods. This preconditioner requires tensor-product bases on quadrilateral or
hexahedral elements. Then, the $p^d \times p^d$ blocks that arise in these
systems can be well-approximated by certain Kronecker products of one
dimensional $p \times p$ matrices. Using a shuffled singular value decomposition
introduced by Van Loan in \cite{VanLoan:1993}, it is possible to compute
decompositions into tensor products of one-dimensional terms that are optimal in
the Frobenius norm. Using these techniques, it is possible to construct an
approximate tensor-product version of the standard block Jacobi preconditioner,
that avoids inverting, or even storing, the large diagonal blocks of the
Jacobian matrix.

In Section \ref{sec:eqns}, we give a very brief description of the discontinuous
Galerkin method for a general system of hyperbolic conservation laws. In Section
\ref{sec:sum-factorization}, we outline the sum-factorization approach, and
describe equivalent Kronecker-product representations. Then, in Section
\ref{sec:tensor-preconds} we develop the approximate Kronecker-product
preconditioners, and provide a new set of algorithms that can be used to
efficiently compute and apply these preconditioners. Finally, in Section
\ref{sec:results}, we apply these preconditioners to several test problems,
including the scalar advection equation, compressible Navier-Stokes equations, 
and the Euler equations of gas dynamics, in two and three spatial dimensions.

%%%%%%%%%%%%%%%%%%%%%%%%%%%%%%%%%%%%%%%%%%%%%%%%%%%%%%%%%%%%%%%%%%%%%%%%%%%%%%%%
% Equations and spatial discretization
%%%%%%%%%%%%%%%%%%%%%%%%%%%%%%%%%%%%%%%%%%%%%%%%%%%%%%%%%%%%%%%%%%%%%%%%%%%%%%%%

\section{Equations and spatial discretization}\label{sec:eqns}
We give a brief overview of the discontinuous Galerkin method for solving a 
hyperbolic conservation law of the form
\begin{equation} \label{eq:cons-law}
  \partial_t u + \nabla\cdot F(u) = 0.
\end{equation}
In order to formulate the method, we first discretize the spatial domain
$\Omega$ by means of a triangulation $\mathcal{T}_h = \{ K_j : \bigcup_j K_j =
\Omega \}$. Common choices for the elements $K_j$ of the triangulation are
simplex and block elements. Given a triangulation $\mathcal{T}_h$, we now
introduce the finite element space $V_h$, given by
\begin{equation}
  V_h = \left\{ v_h : v_h |_{K_j} \in V(K_j) \right\},
\end{equation}
where $V(K_j)$ is a function space local to the element $K_j$. Such functions
admit discontinuities along the element interfaces $\partial K_j$. In the case
of simplex elements, the local function space $V(K_j)$ is taken to be the space
of multivariate polynomials of at most degree $p$, $\mathcal{P}^p(K_j)$.  Of
particular interest to this paper are the block elements, which in
$\mathbb{R}^d$ are defined as the image of the $d$-fold cartesian product of the
interval $[0,1]$ under an isoparametric polynomial transformation map.

By looking for a solution $u_h \in V_h$, multiplying by a test function $v_h \in
V_h$, and integrating by parts over each element, we derive the weak formulation
of 
\eqref{eq:cons-law},
\begin{equation} \label{eq:weak-form}
  \int_{K_j} (\partial_t u_h) v_h~dx
  - \int_{K_j} F(u_h) \cdot \nabla v_h ~dx
  + \int_{\partial K_j} \widehat{F}(u_h^-, u_h^+, n) v_h~dA = 0,
  \qquad\text{for all $K_j \in \mathcal{T}_h$,}
\end{equation}
where $u_h^-$ and $u_h^+$ are the interior and exterior traces (respectively) of
$u_h$ on $\partial K_j$, and $\widehat{F}$ is an appropriately defined
\textit{numerical flux function}. The integrals in \eqref{eq:weak-form} are
approximated using an appropriate quadrature rule, and the resulting system of
ordinary differential equations is termed the semi-discrete system. In this
work, we use quadrature rules that are given by tensor products of
one-dimensional quadratures. Typically, using the method of lines, the time
derivative \eqref{eq:weak-form} is discretized by means of one of many standard
(implicit or explicit) methods for solving ordinary differential equations.

%%%%%%%%%%%%%%%%%%%%%%%%%%%%%%%%%%%%%%%%%%%%%%%%%%%%%%%%%%%%%%%%%%%%%%%%%%%%%%%%
% The sum-factorization approach
%%%%%%%%%%%%%%%%%%%%%%%%%%%%%%%%%%%%%%%%%%%%%%%%%%%%%%%%%%%%%%%%%%%%%%%%%%%%%%%%

\section{The sum-factorization approach} \label{sec:sum-factorization}
In order to numerically represent the solution, we expand the function $u_h$ in
terms of basis functions local to each element. In $\mathbb{R}^d$, the number of
degrees of freedom $n$ per element thus scales as $\mathcal{O}(p^d)$. In this
work, we will make the assumption that the number of quadrature points, denoted
$\mu$, is given by a constant multiple of $n$, and thus also $\mathcal{O}(p^d)$.
We first note that in order to approximate the integrals of the weak form
\eqref{eq:weak-form}, we must evaluate a function $v_h \in V_h$ at each of the
quadrature nodes in a given element $K$. We suppose that 
$\{\Phi_1, \ldots, \Phi_n\}$ is a basis for the local space $V(K)$. Then, for 
any $x \in K$, we can expand $v_h$ in terms of its coefficients $v_j$
\begin{equation}
  \label{eq:gauss-value} v_h(x) = \sum_{j=1}^n v_j \Phi_j(x).
\end{equation}
It is important to note that each computation of $v_h(x_\alpha)$ thus requires
$n$ evaluations of the basis functions. Performing this computation for each
quadrature point therefore requires a total of $\mathcal{O}(p^{2d})$
evaluations. In order to reduce the computational cost of this, and other
operations, we describe the sum-factorization approach, first introduced in 
\cite{Orszag:1980}, and extended to the DG method in \textit{e.g.} 
\cite{Vos:2010}.

\subsection{Tensor-product elements}
The key to the sum-factorization approach are tensor-product elements, where
each element $K$ in the triangulation $\mathcal{T}_h$ is given as the image of
the cartesian product $[0,1]^d$ under a transformation mapping 
{(that is to say, the mesh consists entirely of mapped quadrilateral 
or hexahedral elements)},
and where the
local basis for each element $K$ is given as the product of one-dimensional
basis functions. To be precise, we define the reference element to be the
$d$-dimensional unit cube $\mathcal{R} = [0,1]^d$, and suppose that
$T(\mathcal{R}) = K$, where $T : \mathbb{R}^d \to \mathbb{R}^d$ is an
isoparametric $p$th degree polynomial map. We let $\{ \phi_j(x) \}_{j=0}^p$ be a
basis for $\mathcal{P}^p([0,1])$, the space of polynomials of degree at most $p$
on the unit interval.
{
Then, we define $V(\mathcal{R})$ to be the tensor-product function space, given
as the space of all functions $f(x_1, x_2, \ldots, x_d) = f_1(x_1) f_2(x_2) 
\cdots f_d(x_d)$, written $f = f_1 \otimes f_2 \otimes \cdots \otimes f_d$, 
where $f_k \in \mathcal{P}^p([0,1])$ for each $1\leq k\leq d$. We write the 
tensor-product basis for this space as}
\begin{equation} \label{eq:tensor-basis}
  V(\mathcal{R}) = 
  \bigotimes_{\ell=1}^d \mathcal{P}^p([0,1]) = 
  \mathrm{span} \left\{
   \phi_{i_1} \otimes \phi_{i_2} \otimes \cdots \otimes \phi_{i_d} : 
   0 \leq i_k \leq p
  \right\},
\end{equation}
where $\phi_{i_1} \otimes \phi_{i_2} \otimes \cdots \otimes \phi_{i_d}
(x_1,x_2,\ldots,x_d) = \phi_{i_1}(x_1) \phi_{i_2}(x_2) \cdots \phi_{i_d}(x_d)$.
As a particular example, we consider the one-dimensional nodal basis for
$\mathcal{P}^p([0,1])$, with nodes $B = \{b_1, b_2, \ldots, b_{p+1}\} \subseteq
[0,1]$, and $\phi_j$ is the unique degree-$p$ polynomial such that $\phi_j(b_k)
= \delta_{jk}$. Thus, the coefficients $v_j$ for a function $v_h \in
\mathcal{P}^p([0,1])$ are given by the nodal values $v_j = v_h(b_j)$. The
tensor-product basis defined by \eqref{eq:tensor-basis} consists exactly of the
multivariate polynomials defined by the nodal basis with nodes given by the
$d$-fold cartesian product $B \times \cdots \times B$. In other words, the basis
functions are given by $\Phi_{i_1,i_2,\ldots,i_d}$, which is the unique
multivariate polynomial of degree at most $p$ in each variable, such that
$\Phi_{i_1,i_2,\ldots,i_d} (b_{j_1}, b_{j_2}, \ldots, b_{j_d}) = \delta_{i_1
j_1} \delta_{i_2 j_2} \cdots \delta_{i_d j_d}$.

\begin{figure}[t!]
    \centering
    \begin{subfigure}[t]{1.6in}
        \centering
        \includegraphics[width=1.1in]{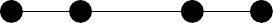}
        %\caption{Cartesian mesh for advection equation}
    \end{subfigure}%
    ~ 
    \begin{subfigure}[t]{1.6in}
        \centering
        \includegraphics[width=1.1in]{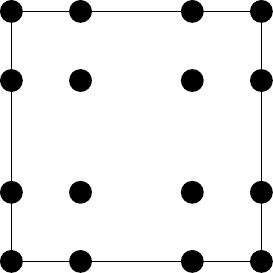}
        %\caption{Unstructured mesh for advection equation}
        %\label{fig:adv-unstructured-mesh}
    \end{subfigure}
    ~ 
    \begin{subfigure}[t]{2.2in}
        \centering
        \includegraphics[width=1.5in]{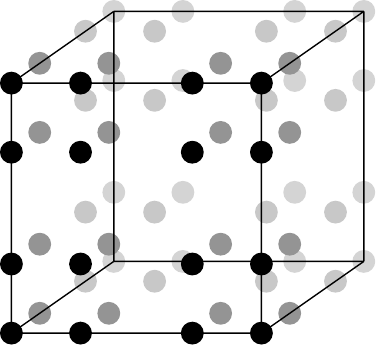}
        %\caption{Unstructured mesh for advection equation}
        %\label{fig:adv-unstructured-mesh}
    \end{subfigure}
    \caption{ 
    Reference elements in 1D, 2D, and 3D with nodes corresponding to 
    $p = 3$. The cartesian-product structure of the nodes gives rise to the 
    tensor-product structure of the corresponding nodal basis functions.}
    \label{fig:ref-elements}
\end{figure}

For a particular element $K = T(\mathcal{R})$, we define the basis for the space
$V(K)$ by means of the transformation map $T$. Given $x \in K$, we write 
$x = T(\xi)$, where $\xi$ denotes the \textit{reference coordinate}. Then, we 
define the basis function $\tilde\Phi_{i_1,i_2,\ldots,i_d}(x)$ by 
$\tilde\Phi_{i_1,i_2,\ldots,i_d}(x) = \Phi_{i_1,i_2,\ldots,i_d}(\xi)$. 
Similarly, it will often be convenient to identify a given function 
$\tilde{v} \in V(K)$ with the function $v \in V(\mathcal{R})$ obtained by 
$v = \tilde v \circ T^{-1}$.

Given this choice of basis, we also define the quadrature nodes on the
$d$-dimensional unit cube to be the $d$-fold cartesian product of given
one-dimensional quadrature nodes, whose weights are the corresponding products
of the one-dimensional weights. Equipped with these choices, we return to the
calculation of the quantity \eqref{eq:gauss-value}. For the sake of
concreteness, we consider the case where $d=3$, for which the calculation above
should naively require $\mathcal{O}(p^6)$ evaluations. We suppose that the
one-dimensional quadrature points are given as $x_1, x_2, \ldots, x_\mu$, and
hence we can write the three-dimensional quadrature points as $\bm{x}_{\alpha,
\beta, \gamma} = (x_\alpha, x_\beta, x_\gamma)$, for all $1 \leq \alpha, \beta,
\gamma \leq \mu$. We then factor the summation in \eqref{eq:gauss-value} to
obtain
\begin{equation} \label{eq:quad-sum-factorization}
\begin{aligned}
  v_h(\bm{x}_{\alpha,\beta,\gamma})
    &= \sum_{i,j,k=1}^{p+1} v_{ijk} \Phi_{ijk}(x_\alpha, x_\beta, x_\gamma) \\
    &= \sum_{k=1}^{p+1} \sum_{j=1}^{p+1} \sum_{i=1}^{p+1}
        v_{ijk} \phi_i(x_\alpha) \phi_j(x_\beta) \phi_k(x_\gamma) \\
    &= \sum_{k=1}^{p+1} \phi_k(x_\gamma) \sum_{j=1}^{p+1} \phi_j(x_\beta)
       \sum_{i=1}^{p+1} v_{ijk} \phi_i (x_\alpha).
\end{aligned}
\end{equation}
We notice that the index of each summation ranges over $p+1$ values, and there
are three free indices in each sum. Thus, the total number of operations
required to evaluate a function $v_h$ at each of the quadrature points is
$\mathcal{O}(p^4)$. For general dimension $d$, this computation requires
$\mathcal{O}(p^{d+1})$ basis function evaluations. Thus, the computational work
\textit{per degree of freedom} is linear in the degree $p$ of polynomial basis,
in contrast to the original estimate of $\mathcal{O}({p^d})$ computational work
per degree of freedom, which is exponential in spatial dimension $d$. In a
similar fashion, the tensor-product structure of this function space can be
exploited in order to compute the integrals in equation \eqref{eq:weak-form} in
linear time per degree of freedom.

\subsection{Kronecker-product structure} \label{sec:kron-structure}
The sum-factorization procedure shown in \eqref{eq:quad-sum-factorization} can 
be described simply and elegantly as a linear-algebraic Kronecker product 
structure. We recall that the Kronecker product of matrix $A^{k \times \ell}$ 
and $B^{m \times n}$ (whose dimensions are indicated by the superscripts), 
is the $km \times \ell n$ matrix $C$ defined by
\begin{equation}
  A \otimes B = C = \left(
    \begin{array}{cccc}
      a_{11} B & a_{12} B & \cdots & a_{1\ell} B \\
      a_{21} B & a_{22} B & \cdots & a_{2\ell} B \\
      \vdots   & \vdots   & \ddots & \vdots      \\
      a_{k1} B & a_{k2} B & \cdots & a_{k\ell} B
    \end{array}
  \right).
\end{equation}
The Kronecker product has many desirable and useful properties, enumerated 
in Van Loan's exposition \cite{VanLoan:2000}.

We can define the one-dimensional \textit{Gauss point evaluation matrix} 
as the $\mu \times (p+1)$ Vandermonde-type matrix obtained by evaluating each 
of the one-dimensional basis functions at all of the quadrature points,
\begin{equation}
  G_{\alpha j} = \phi_j(x_\alpha),
\end{equation}
and, in a similar fashion, it is also convenient to define the 
one-dimensional differentiation matrix, $D$, whose entries are given by
$ D_{\alpha j} = \phi_j'(x_\alpha)$. We now describe the Kronecker-product 
structure of a general $d$-dimensional DG method, for arbitrary $d$. 
Let $I = (i_1, i_2, \ldots, i_d)$ and 
$A = (\alpha_1, \alpha_2, \ldots, \alpha_d)$ be multi-indices of length $d$.
Then, we can define a vector $v$ of length $(p+1)^d$ whose entries are given by 
concatenating the entries of the $d$th-order tensor $v_I$. Thus, we obtain the 
values of $v_h$ evaluated at the quadrature points $\bm{x}_A$ by computing
the $d$-fold Kronecker product
\begin{equation}
  v_h(\bm{x}_{A}) = (G \otimes G \otimes \cdots \otimes G) v.
\end{equation}
This Kronecker-product representation is computationally equivalent to the
sum-factorized version from the preceding section. Indeed, many of the
operations needed for the computation of the discontinuous Galerkin method are
amenable to being written in the form of Kronecker products.

For instance, it is often useful to approximate quantities of the form
\begin{equation}
  \int_{K} f(x) v_h(x)~dx,
\end{equation}
where $f$ is an arbitrary function whose value is known at the appropriate 
quadrature nodes. This requires approximating the integrals
\begin{equation} \label{eq:int-f-phi}
  \int_{K} f(x) \tilde\Phi_{I}(x)~dx
\end{equation}
for each of the basis functions $\tilde\Phi_{I}$. We consider the element $K$
to be the image under the isoparametric transformation map $T$ of the reference
element $\mathcal{R} = [0,1]^d$. In this notation, for all $x \in K$, $x =
T(\xi)$, where $\xi \in \mathcal{R}$. We then write \eqref{eq:int-f-phi} as an
integral over the reference element,
\begin{equation} \label{eq:int-f-phi-ref}
  \int_{K} f(x) \tilde\Phi_{I}(x)~dx = 
  \int_{\mathcal{R}} f(T(\xi)) \Phi_{I}(\xi) | \det(T'(\xi))|~d\xi.
\end{equation}
To this end, we define a diagonal
weight matrix $W$ by whose entries along the diagonal are given by $w_\alpha$,
where $w_\alpha$ is the quadrature weight associated with the point $x_\alpha$.
Additionally, we define the $\mu^d \times \mu^d$ diagonal matrix $J_T$ whose 
entries are equal to the Jacobian determinant of the isoparametric mapping 
$\det(T'(\xi))$ at each of the quadrature points. Then, the $(p+1)^d$ integrals
of the form \eqref{eq:int-f-phi-ref} can be found as the entries of the vector
\begin{equation}
  \left(G^T W\right) \otimes \left(G^T W\right) \otimes
  \cdots \otimes \left(G^T W\right)
  J_T f(\bm{x}_{A}).
\end{equation}
In a similar fashion, the computation of all of the quantities needed to 
formulate a DG method can be written in Kronecker form. In Table 
\ref{tab:dg-ops}, we summarize the Kronecker-product formulation of 
several other important operations needed for the DG method, for 
the special cases of $d=2$ and $d=3$.

\begin{table}[t]
\caption{Kronecker-product form of DG operations}
\label{tab:dg-ops}
{\def\arraystretch{1.75}\tabcolsep=3pt
\begin{tabular}{p{2.5in}ll}
\toprule
Operation & 2D & 3D \\
\midrule
\raggedright
Evaluate solution at quadrature points
  & $\begin{aligned} \left(G \otimes G\right) u \end{aligned} $
  & $\begin{aligned} \left(G \otimes G \otimes G\right) u \end{aligned} $ \\
\raggedright
Integrate function $f$ (known at quadrature points) against test functions
  & $\begin{aligned} \left(G^T W \otimes G^T W \right) J_T f \end{aligned}$
  & $\begin{aligned} \left(G^T W \otimes G^T W \otimes G^T W \right) J_T f
     \end{aligned}$ \\
\raggedright
Integrate function $f = (f_1,\ldots,f_d)$ against gradient of test functions
  & $\begin{aligned}[t]
    \left(G^T W \otimes D^T W \right) J_T f_1\\
    \left(D^T W \otimes G^T W \right) J_T f_2
    \end{aligned}$
  & $ \begin{aligned}[t]
    \left(G^TW \otimes G^T W \otimes D^T W \right) J_T f_1\\
    \left(G^TW \otimes D^T W \otimes G^T W \right) J_T f_2\\
    \left(D^TW \otimes G^T W \otimes G^T W \right) J_T f_3
   \end{aligned} $\\
\bottomrule
\end{tabular}}
\end{table}

\subsection{Explicit time integration} \label{sec:explicit}
It is important to note that all of these operations have a computational 
complexity of at most $\mathcal{O}(p^{d+1})$. In other words, the cost of these 
operations scales linearly in $p$ per degree of freedom. We now return to the 
semi-discrete system of equations \eqref{eq:weak-form}, which we rewrite as
\begin{equation} \label{eq:mass-u-r}
  M (\partial_t u_h) = r,
\end{equation}
where $r$ is the quadrature approximation of the second two integrals on the 
left-hand side of \eqref{eq:weak-form}. Then, using the operations described 
above, it is possible to compute all the integrals required to form $r$. What is
required in order to integrate this semi-discrete equation explicitly in time 
is to invert the mass matrix $M$. We recall that, with a tensor-product 
basis, we can compute the element-wise mass matrix (on, \textit{e.g.}, element
$K_j$) as
\begin{equation} \label{eq:mass-kron}
  M_{j} = \left((G^T W) \otimes \cdots \otimes (G^T W)\right) J_T
  \left(G \otimes \cdots \otimes G\right)
\end{equation}
where $G$ is the Gauss point evaluation matrix defined above, $W$ is the 
diagonal matrix with the one-dimensional quadrature weights on the diagonal, 
and $J_T$ is the $(\mu^d) \times (\mu^d)$ diagonal matrix whose entries are 
equal to the absolute Jacobian determinant of the element transformation map, 
evaluated at each of the quadrature points.

One strategy, proposed in \cite{NME:NME5137}, is to use the same number of 
quadrature points as DG nodes, such that $\mu = p+1$. In that case, all of 
the matrices appearing in \eqref{eq:mass-kron} are square, and we can compute
\begin{equation} \label{eq:minv-kron}
  M_j^{-1} = \left(G^{-1} \otimes \cdots \otimes G^{-1} \right)
           J_T^{-1}
           \left((G^T W)^{-1} \otimes \cdots \otimes (G^T W)^{-1}\right).
\end{equation}
Since $G$ and $G^T W$ are both $(p+1)\times(p+1)$ matrices, these operations 
can be performed in $\mathcal{O}(p^3)$ time. Additionally, $J_T$ is a
$(p+1)^d \times (p+1)^d$ diagonal matrix, and thus can be inverted in 
$(p+1)^d$ operations. (On a practical note, in this case we would avoid 
explicitly forming the inverse matrices $G^{-1}$ and $(G^TW)^{-1}$, and would 
instead opt to form their $LU$ factorizations). Thus, the linear system 
\eqref{eq:mass-u-r} can be solved in the same complexity as multiplying by the 
expression on the right-hand side of \eqref{eq:minv-kron}, \textit{i.e.} 
$\mathcal{O}(p^{d+1})$.

A serious drawback to this approach is that using the same number of quadrature
points as DG nodes does not, in general, allow for exact integration of the 
quantity
\begin{equation} \label{eq:mass-int}
  \int_K u_h v_h~dx, \qquad u_h, v_h \in V_h,
\end{equation}
because of the use of isoparametric elements, where the Jacobian determinant of
the transformation mapping may itself be a high-degree polynomial. In order to
address this issue, we introduce a new strategy for solving the system
\eqref{eq:mass-u-r}. We first note that the global mass matrix has a natural
element-wise block-diagonal structure, where, furthermore, each block $M_j$ is a
symmetric positive-definite matrix. Thus, we can solve this system of equations
element-by-element, using the preconditioned conjugate gradient (PCG) method. As
a preconditioner, we use the under-integration method described above.

Thus, each iteration in the PCG solver requires a multiplication by the exact
mass matrix, and a linear solve using the approximate, under-integrated mass
matrix, given by equation \eqref{eq:minv-kron}, for the purposes of
preconditioning. Both of these operations are performed in
$\mathcal{O}(p^{d+1})$ time by exploiting their tensor-product structure.
This has the consequence that if the element transformation mapping is bilinear
(and so the corresponding element is straight-sided) then its Jacobian
determinant is linear, and with the appropriate choice of quadrature points, the
integral \eqref{eq:mass-int} can be computed exactly, and the PCG method will
converge within one iteration. In practice, we observe that the number of PCG
iterations required to converge is very small even on curved, isoparametric
meshes, and does not grow with $p$.

The techniques described above are sufficient to implement an explicit
discontinuous Galerkin method with tensor-product elements, requiring
$\mathcal{O}(p^{d+1})$ operations per time step.
The main restriction to using such explicit methods with very high polynomial
degree $p$ is the restrictive CFL condition. It has been shown that the rate of
growth of the spectral radius of the semi-discrete DG operator is bounded above
by $(p+1)(p+2)/h$ \cite{Warburton2008,gottlieb1991cfl}, and well-approximated by
$(p+1)^{1.78}$ \cite{Krivodonova:2013}. This requires that the time step satisfy
approximately approximately $\Delta t\leq Ch/(p+1)^{1.78}$, which can prove to
be prohibitively expensive as the number of time steps needed increases.
For this reason, we are interested in applying some of
the same tensor-product techniques to efficiently integrate in time implicitly,
and thus avoid the overly-restrictive CFL condition.

\subsection{Implicit time integration}
Instead of an explicit time integration method, we now consider an implicit 
schemes such as backward differentiation formulas (BDF) or diagonally-implicit 
Runge-Kutta (DIRK) methods. The main advantage of such methods is that they 
remain stable for larger time steps, even in the presence of highly anisotropic 
elements. Additionally, these methods avoid the restrictive $p$-dependent 
explicit stability condition mentioned above. Such implicit methods require the 
solution of systems of the form
\begin{equation}
  M u_h - \Delta t f(u_h) = r,
\end{equation}
which, when solved by means of Newton's method, give rise to linear systems of 
the form
\begin{equation}
  (M - \Delta t J) x = b,
\end{equation}
where the matrix $J$ is the Jacobian of the potentially non-linear function $f$.

The most immediate challenge towards efficiently implementing an implicit method
for high polynomial degree on tensor product elements is forming the Jacobian
matrix. In general, all the degrees of freedom within one element are coupled,
and thus the diagonal blocks of the Jacobian matrix corresponding to a single
element are dense $(p+1)^d \times (p+1)^d$ matrices. Therefore, it is impossible
to explicitly form this matrix in less than $\mathcal{O}(p^{2d})$ time. To
circumvent this, using the techniques described in the preceding sections, it is
possible to solve the linear systems arising from implicit time integration by
means of an iterative method such as GMRES \cite{Persson:2008}. Each iteration
requires performing a matrix-vector multiplication by the mass matrix and the
Jacobian matrix. If we avoid explicitly forming these matrices, then the
multiplications can be performed in $\mathcal{O}(p^{d+1})$ time, using methods
similar to those described in the preceding section.

\subsubsection{Matrix-free tensor-product Jacobians}\label{sec:matrix-free}
In order to efficiently implement the implicit method described above, we first 
apply the sum-factorization technique to efficiently evaluate the matrix-vector
product
\begin{equation}
  (M - \Delta t J) v
\end{equation}
for a given vector $v$. As described in Section \ref{sec:kron-structure}, the
mass matrix is a block diagonal matrix whose $j$th block can be written in
the Kronecker form given by equation \eqref{eq:mass-kron}.
Multiplying by the Kronecker products can be performed with
$\mathcal{O}(p^{d+1})$ operations, and multiplying by $J_T$ requires exactly
$(p+1)^d$ operations. Thus, the product $M x$ requires $\mathcal{O}(p^{d+1})$
operations.

We now describe our algorithm for also computing the product $J x$ with the same
complexity. For simplicity of presentation, we will take $d=2$, but the
algorithm is immediately generalizable to arbitrary dimension $d$. We first
consider an element-wise blocking of the matrix $J$. Each block is a $(p+1)^2
\times (p+1)^2$ matrix, with blocks along the diagonal corresponding to each
element in the triangulation, and blocks off the diagonal corresponding to the
coupling between neighboring elements through their common face. We write $r =
f(u)$, and restrict our attention to one element. We define the indices $1 \leq
i,j,k,\ell \leq p+1$, such that the entries of the diagonal block of the
Jacobian can be written as
\begin{equation}
  \frac{\partial r_{ij}}{\partial u_{k\ell}}.
\end{equation}
We define $r_{ij}$ by
\begin{equation}
r_{ij} = \int_{K} F(u) \cdot \nabla \tilde\Phi_{ij} ~dx
  - \int_{\partial K} \widehat{F}(u^-, u^+, n) \tilde\Phi_{ij}~dA
\end{equation}
which we evaluate using the quadrature rule
\begin{multline}
  r_{ij} = \sum_{\alpha = 1}^\mu \sum_{\beta = 1}^\mu w_\alpha w_\beta
           F(u(x_\alpha, x_\beta)) \cdot \nabla 
           \left(\phi_i(x_\alpha) \phi_j(x_\beta)\right) \\
           - \sum_{e \in \partial K}
           \sum_{\alpha = 1}^\mu  w_\alpha \widehat{F}(u^-(x^e_\alpha, y^e_\alpha),
           u^+(x^e_\alpha, y^e_\alpha), n(x^e_\alpha, y^e_\alpha))
           \phi_i(x^e_\alpha) \phi_j(y^e_\alpha),
\end{multline}
where the notation $(x^e_\alpha, y^e_\alpha)$ represents the coordinates of the
$\alpha$th quadrature node alone the face $e$ of $\partial K$. We also recall
that the function $u$ is evaluated by expanding in terms of the local basis
functions, \textit{e.g.}
\begin{equation}
  u(x_\alpha, x_\beta) = \sum_{k=1}^{p+1} \sum_{\ell=1}^{p+1}
    u_{k\ell} \phi_k(x_\alpha) \phi_\ell(x_\beta),
\end{equation}
which can be evaluated efficiently as
\begin{equation}
  \left(G \otimes G\right) u.
\end{equation}
Thus, the entries of the Jacobian can be written as
\begin{multline} \label{eq:jac-entries}
  \frac{\partial r_{ij}}{\partial u_{k\ell}} = 
       \sum_{\alpha = 1}^\mu \sum_{\beta = 1}^\mu w_\alpha w_\beta
            \phi_k(x_\alpha) \phi_\ell(x_\beta)
           \frac{\partial F}{\partial u}(u(x_\alpha, x_\beta)) \cdot \nabla 
           \left(\phi_i(x_\alpha) \phi_j(x_\beta)\right) \\
           - \sum_{e \in \partial K}
           \sum_{\alpha = 1}^\mu w_\alpha\phi_k(x^e_\alpha)\phi_\ell(y^e_\alpha)
           \frac{\partial\widehat{F}}{\partial u^-}
           (u^-(x^e_\alpha, y^e_\alpha),
           u^+(x^e_\alpha, y^e_\alpha), n(x^e_\alpha, y^e_\alpha))
           \phi_i(x^e_\alpha) \phi_j(y^e_\alpha).
\end{multline}
Since there are $(p+1)^{2d}$ such entries, we avoid explicitly computing the 
entries of this matrix, and instead describe how to compute the matrix-vector 
product $Jx$. As a pre-computation step, we compute the flux derivatives 
$\partial F / \partial u$ and numerical flux derivatives
$\partial \widehat{F} / \partial u^-$ at each of the quadrature nodes 
$(x_\alpha, x_\beta)$ in the element $K$. For simplicity, we introduce the 
notation
\begin{equation}
\frac{\partial F}{\partial u} (x_\alpha, x_\beta)
  =\frac{\partial F}{\partial u}(u(x_\alpha, x_\beta)), \qquad
\frac{\partial\widehat{F}}{\partial u^-} (x^e_\alpha, y^e_\alpha)
  = \frac{\partial\widehat{F}}{\partial u^-}
    (u^-(x^e_\alpha, y^e_\alpha),
    u^+(x^e_\alpha, y^e_\alpha), n(x^e_\alpha, y^e_\alpha)).
\end{equation}
Then, applying the sum-factorization technique, the terms of the product of $Jv$
for a given vector $v$ corresponding to the diagonal block takes the form
\begin{align}
\begin{split}
  \left(\frac{\partial r_{ij}}{\partial u_{k\ell}}\right) v_{k\ell}
   &= \sum_{k=1}^{p+1} \sum_{\ell=1}^{p+1} 
      \sum_{\alpha = 1}^\mu \sum_{\beta = 1}^\mu
          w_\alpha w_\beta \phi_k(x_\alpha) \phi_\ell(x_\beta)
           \frac{\partial F}{\partial u}(x_\alpha, x_\beta) \cdot \nabla 
           \left(\phi_i(x_\alpha) \phi_j(x_\beta)\right) v_{k \ell} \\
  & \qquad - \sum_{k=1}^{p+1} \sum_{\ell=1}^{p+1}  \sum_{e \in \partial K}
           \sum_{\alpha = 1}^\mu w_\alpha \phi_k(x^e_\alpha)
           \phi_\ell(y^e_\alpha)
           \frac{\partial\widehat{F}}{\partial u^-}
           (x^e_\alpha, y^e_\alpha)
           \phi_i(x^e_\alpha) \phi_j(y^e_\alpha) v_{k \ell}
\end{split} \\
\begin{split}
&=
  \sum_{\alpha=1}^\mu w_\alpha  \phi_i'(x_\alpha)
    \sum_{\beta=1}^\mu w_\beta 
    \frac{\partial F_1}{\partial u} (x_\alpha, x_\beta) \phi_j(x_\beta)
    \sum_{\ell=1}^{p+1} \phi_\ell(x_\beta) 
    \sum_{k=1}^{p+1} \phi_k(x_\alpha) v_{k \ell} \\
& \qquad + \sum_{\alpha=1}^\mu w_\alpha  \phi_i(x_\alpha)
    \sum_{\beta=1}^\mu w_\beta
    \frac{\partial F_2}{\partial u} (x_\alpha, x_\beta) \phi_j'(x_\beta)
    \sum_{\ell=1}^{p+1} \phi_\ell(x_\beta)
    \sum_{k=1}^{p+1} \phi_k(x_\alpha) v_{k \ell} \\
& \qquad + \sum_{e \in \partial K}
           \sum_{\alpha = 1}^\mu w_\alpha 
           \frac{\partial\widehat{F}}{\partial u^-}
           (x^e_\alpha, y^e_\alpha)
           \phi_i(x^e_\alpha) \phi_j(y^e_\alpha) 
           \sum_{\ell=1}^{p+1} \phi_\ell(y^e_\alpha)
           \sum_{k=1}^{p+1} \phi_k(x^e_\alpha) v_{k \ell},
\end{split}
\end{align}
where $F_1$ and $F_2$ are the $x$ and $y$ components of the flux function $F$,
respectively. We notice that in each of the above summations, there are at most
two free indices, and therefore each sum can be computed in $\mathcal{O}(p^3)$
time, achieving linear time in $p$ per degree of freedom. The terms of the
product corresponding to the off-diagonal blocks have a similar form to the face
integral in the above equations, and can similarly be computed in
$\mathcal{O}(p^3)$ time.

To summarize, we describe the algorithm for computing the matrix-products 
of the form $Jx$ in Algorithm \ref{alg:matrix-free-matvec}.

\algblockdefx{Precomp}{EndPrecomp}
  {\textbf{Pre-computation:}}
  {EndPrecomp}
\algnotext{EndPrecomp}
\begin{algorithm}[H]
\caption{ Matrix-free computation of $Jv$ in 2D and 3D}
\label{alg:matrix-free-matvec}
\begin{algorithmic}[1] 
  \Precomp
  \State Evaluate the solution at quadrature points:
  \Statex \hspace{0.8cm} (in 2D, compute 
     $ \left( G \otimes G \right) u $, and in 3D, compute
     $ \left( G \otimes G \otimes G \right) u $)
     \Comment{Complexity: $\mathcal{O}(p^{d+1})$}
  \State Evaluate the flux Jacobians
     $\frac{\partial F}{\partial u}$ and
     $\frac{\partial \widehat{F}}{\partial u^-}$ at quadrature points
     \Comment{Complexity: \makebox[0pt][l]{$\mathcal{O}(p^d)$}%
\phantom{$\mathcal{O}(p^{d+1})$}}
   \EndPrecomp
  \State Compute the matrix-vector product using the sum-factorized form
    \Comment{Complexity: $\mathcal{O}(p^{d+1})$}
\end{algorithmic}
\end{algorithm}

\noindent
The first two operations can be performed as a pre-computation step, and only 
the third step need be repeated when successively multiplying the same Jacobian 
matrix by different vectors (as in the case of an iterative linear solver).

%%%%%%%%%%%%%%%%%%%%%%%%%%%%%%%%%%%%%%%%%%%%%%%%%%%%%%%%%%%%%%%%%%%%%%%%%%%%%%%%
% Tensor-product preconditioners
%%%%%%%%%%%%%%%%%%%%%%%%%%%%%%%%%%%%%%%%%%%%%%%%%%%%%%%%%%%%%%%%%%%%%%%%%%%%%%%%

\section{Tensor-product preconditioners} \label{sec:tensor-preconds}
One of the main challenges in successfully applying such a matrix-free method is
preconditioning \cite{Saad:2003iterative}. Common preconditioners typically used
for implicit DG methods include block Jacobi, block Gauss-Seidel, and block ILU
preconditioners \cite{Persson:2008}. Computing these preconditioners first
requires forming the matrix, and additionally requires the inversion of certain
blocks. Typically, this would incur a cost of $\mathcal{O}(p^{3d})$, which
quickly grows to be prohibitive as we take $p$ to be large. To remedy this
issue, we develop a preconditioner for two and three spatial dimensions that
takes a similar Kronecker product form to those seen in the previous section.

We draw inspiration from the tensor-product structure often seen in
finite-difference and spectral approximations to, \textit{e.g.} the Laplacian
operator on a $n^d$ cartesian grid, which can be written in one, two, and three
spatial dimensions, respectively, as
\begin{equation} \label{eq:poi-tensor}
  L_{\rm 1D} = T_n,\qquad 
  L_{\rm 2D} = I \otimes T_n + T_n \otimes I, \qquad 
  L_{\rm 3D} = I \otimes I \otimes T_n + I \otimes T_n \otimes I
             + T_n \otimes I \otimes I,
\end{equation}
where $T_n$ is the standard one-dimensional approximation to the Laplacian. 
Given a general conservation law of the form \eqref{eq:cons-law}, the flux 
function $F$ is not required to possess any particular structure, and thus the 
DG discretization of such a function will not be exactly expressible in a 
similar tensor-product form. That being said, many of the key operations in DG, 
listed in Table \ref{tab:dg-ops}, are expressible in a similar form. 
Therefore, in order to precondition the implicit systems of the form
\begin{equation}
  (M - \Delta t J) u = b,
\end{equation}
we look for tensor-product approximations to the diagonal blocks $A$ of the
matrix $M - \Delta t J$. 
{
Specifically, we are interested in finding preconditioners $P$ of the form
\begin{alignat}{2} \label{eq:2d-kron-sum}
  A&\approx P=\sum_{j=1}^{r} A_j \otimes B_j &&\text{in 2D}, \\
  \label{eq:3d-kron-sum}
  A&\approx P=\sum_{j=1}^{r} A_j\otimes B_j\otimes C_j \qquad&&\text{in 3D},
\end{alignat}
for a fixed number of terms $r$, where each of the matrices $A_j, B_j,$ and
$C_j$ are of size $(p+1) \times (p+1)$.}
Given $r$, it is possible to find the best possible approximation (in the
Frobenius norm) of the form \eqref{eq:2d-kron-sum} to an arbitrary given matrix
by means of the Kronecker-product singular value decomposition (KSVD).

\subsection{Kronecker-product singular value decomposition}
In  \cite{VanLoan:1993}, Van Loan posed the \textit{nearest Kronecker product
problem} (NKP): given a matrix $A \in \mathbb{R}^{m \times n}$ (with $m = m_1
m_2$ and $n=n_1n_2$), and given a fixed number $r$, find matrices $A_j \in
\mathbb{R}^{m_1 \times n_1}, B_j \in \mathbb{R}^{m_2 \times n_2}$ that minimize
the Frobenius norm
\begin{equation} \label{eq:NKP}
  \left\| A - \sum_{j=1}^r A_j \otimes B_j \right\|_F.
\end{equation}
The solution to the NKP given by Van Loan is as follows. We first consider the 
``blocking'' of $A$:
\begin{equation}
  A = \left(\begin{array}{cccc}
    A_{11}    & A_{12}    & \cdots & A_{1,n_1} \\
    A_{21}    & A_{22}    & \cdots & A_{2,n_1} \\
    \vdots    & \vdots    & \ddots & \vdots    \\
    A_{m_1,1} & A_{m_1,2} & \cdots & A_{m_1,n_1}
  \end{array}\right),
\end{equation}
where each block is a $m_2 \times n_2$ matrix. We then define a rearranged (or
\textit{shuffled}) version $\widetilde{A}$ of the matrix $A$, which is a $m_1n_1
\times m_2n_2$ matrix given by
\begin{equation}
  \widetilde{A} = \left( \begin{array}{c}
    \widetilde{A}_1 \\
    \widetilde{A}_2 \\
    \vdots          \\
    \widetilde{A}_{n_1}
  \end{array}\right),
  \text{ where $\widetilde{A}_j$ is a block of rows given by } 
  \widetilde{A}_j = \left( \begin{array}{c}
    \Vec(A_{1j})^T    \\
    \Vec(A_{2j})^T    \\
    \vdots            \\
    \Vec(A_{m_1,j})^T
  \end{array}\right),
\end{equation}
where the $\Vec$ operator is defined so that $\Vec(A_{ij})$ is the column vector
of length $m_2n_2$ obtained by ``stacking'' the columns of $A_{ij}$. This 
rearranged matrix has the property that, given matrices $A_j, B_j$
\begin{equation}
  \left\| A - \sum_{j=1}^r A_j \otimes B_j \right\|_F 
  = \left\| \widetilde{A} - \sum_{j=1}^r \Vec(A_j) \Vec(B_j)^T \right\|_F,
\end{equation}
and therefore the NKP problem \eqref{eq:NKP} has been reduced the finding the 
closest rank-$r$ approximation to $\widetilde{A}$. This approximation can be 
found by computing the singular value decomposition (SVD) of $\widetilde{A}$, 
\begin{equation}
  \widetilde{A} = U \Sigma V^T,
\end{equation}
then the solution to \eqref{eq:NKP} is given by reshaping the columns of $U$ 
and $V$, such that
\begin{equation}
  \Vec(A_j) = \sqrt{\sigma_j} U_j, \qquad \Vec(B_j) = \sqrt{\sigma_j} V_j.
\end{equation}
This construction is referred to as the \textit{Kronecker product SVD} (KSVD).

\subsubsection{Efficient computation of the KSVD}
In general, computing the singular value decomposition of a matrix is an
expensive process, with cubic complexity. However, if the number $r$ of desired
terms in the summation \eqref{eq:2d-kron-sum} is much smaller than the rank
$\tilde{r}$ of the matrix $\widetilde{A}$, then it is possible to
well-approximate the largest singular values and associated left and right
singular vectors by means of a Lanczos algorithm \cite{Golub:1981}. This
algorithm has the additional advantage that an explicit representation of the
matrix $\widetilde{A}$ is not required, rather only the ability to multiply
vectors by the shuffled matrices $\widetilde{A}$ and $\widetilde{A}^T$. In this
section, we follow the presentation from \cite{VanLoan:1993}. The Lanczos SVD
procedure is described in Algorithm \ref{alg:lanczos}.

\begin{algorithm}
\caption{Lanczos SVD}
\label{alg:lanczos}
\begin{algorithmic}[1]
  \State $v_0 \gets $ random vector with $\|v_0\|_2 = 1$
  \State $p_0 \gets v_0, \beta_0 \gets 1, u_0 \gets 0$
  \For{$j=0$ to $J$ (maximum number of iterations)}
    \State $v_{j} \gets p_j/\beta_j$
    \State $r_{j} \gets \widetilde{A}v_{j} - \beta_{j}u_{j}$
      \label{alg:step-mv}
    \State Orthogonalize $r_j$. \label{alg:lanczos:orth1}
    \State $\alpha_{j} \gets \|r_{j}\|_2$
    \State $u_{j+1} \gets r_j/\alpha_j$
    \State $p_{j+1} \gets \widetilde{A}^T u_j - \alpha_{j}v_{j}$
    \State Orthogonalize $p_{j+1}$. \label{alg:lanczos:orth2}
      \label{alg:step-mvT}
    \State $\beta_{j+1} \gets \|p_{j+1}\|_2$
    \If{$|\beta_{j+1}| < \rm{tolerance}$}
      \State \textbf{break}
    \EndIf
  \EndFor
  \State $U \gets ( u_1,  u_2, \ldots, u_{j+1} )$
  \State $V \gets ( v_0,  v_1, \ldots, v_j )$
  \State Construct bidiagonal matrix $B$, with diagonal $\alpha_0, \ldots, 
         \alpha_j$, and superdiagonal $\beta_1, \ldots, \beta_j$.
  \State Compute $r$ largest singular values $\sigma_k$ (and corresponding
         left and right singular vectors, $u_k'$, $v_k'$) of $B$
  \State Singular values of $\widetilde{A}$ are $\sigma_j$, singular vectors
         are $U u_k'$ and $V v_k'$.
\end{algorithmic}
\end{algorithm}

We remark that there are many variations on the orthogonalization procedure 
referred to in lines \ref{alg:lanczos:orth1} and \ref{alg:lanczos:orth2} of 
Algorithm \ref{alg:lanczos}, including partial or full orthogonalization. In 
this work, we perform full orthogonalization of the vectors $u_j$ and $v_j$ 
at each iteration of the Lanczos algorithm.

As mentioned previously, one of the key advantages of the Lanczos algorithm is
that an explicit representation of the matrix $\widetilde{A}$ can be foregone,
since only matrix-vector products of the form $\widetilde{A}x$ and
$\widetilde{A}^T x$ are required. As described in \cite{VanLoan:1993}, we can
compute these matrix-vector products according to Algorithms \ref{alg:matvec}
and \ref{alg:matvecT}. Taking advantage of the specific tensor-product form of
the matrix $A$, and using techniques similar to those used for the matrix-free
Jacobian evaluation from Section \ref{sec:matrix-free}, it is possible to
efficiently compute the matrix-vector products. Specialized kernels are required
for two and three spatial dimensions, and the details of this process are
described in the following sections.

\begin{algorithm}
\caption{Compute $u = \widetilde{A}v$}
\label{alg:matvec}
\begin{algorithmic}[1]
   \State $u \gets 0$
   \For{$i=1$ to $n_1$}
     \State $\mathrm{rows} \gets (i-1)m_1 + 1, \ldots, im_1$
     \For{$j=1$ to $n_2$}
       \State Define $Z \in \mathbb{R}^{m_2 \times m_1}$ by
              $\Vec(Z) = A(\ :\ ,(i-1)n_2 + j)$
       \State
        $u(\mathrm{rows}) \gets u(\mathrm{rows}) + Z^T v((j-1)m_2 + 1{:}jm_2)$
     \EndFor
   \EndFor
\end{algorithmic}
\end{algorithm}
\begin{algorithm}
\caption{Compute $u = \widetilde{A}^Tv$}
\label{alg:matvecT}
\begin{algorithmic}[1]
   \State $u \gets 0$
   \For{$i=1$ to $n_2$}
     \State $\mathrm{rows} \gets (i-1)m_2 + 1, \ldots, im_2$
     \For{$j=1$ to $n_1$}
       \State Define $Z \in \mathbb{R}^{m_2 \times m_1}$ by
              $\mathrm{vec}(Z) = A(\ :\ ,(j-1)n_2 + i)$
        \State
          $u(\mathrm{rows}) \gets u(\mathrm{rows}) + Z v((j-1)m_1 + 1{:}jm_1)$
     \EndFor
   \EndFor
\end{algorithmic}
\end{algorithm}

\subsection{Two spatial dimensions}
Having shown that, given the number of terms $r$ in the sum, it is possible to 
find the best approximation of the form \eqref{eq:2d-kron-sum}, we now address 
the issue of solving linear systems of equations with such a matrix. In the case
that $r=1$, we have $P = A_1 \otimes B_1$, and it is clear that 
$P^{-1} = A_1^{-1} \otimes B_1^{-1}$, and thus the $(p+1)^2 \times (p+1)^2$ 
problem is reduced to two problems of size $(p+1) \times (p+1)$. Our experience 
has shown that $r=1$ is not sufficient to accurately approximate the Jacobian 
matrix, and the resulting preconditioners are not very effective. For this 
reason, we choose $r=2$, and obtain a linear system of the form
\begin{equation}
  Px = (A_1 \otimes B_1 + A_2 \otimes B_2)x = b.
\end{equation}
Because of the additional term in this sum, it is not possible to invert this
matrix factor-wise. Instead, we follow the matrix diagonalization technique
described in \cite{Shen2011, Lynch:1964direct}. We multiply the system of
equations on the left by $(A_2^{-1} \otimes B_1^{-1})$ to obtain
\begin{equation}
  (A_2^{-1}A_1 \otimes I + I \otimes B_1^{-1}B_2) x = 
  (A_2^{-1} \otimes B_1^{-1})b.
\end{equation}
We let $C_1 = A_2^{-1}A_1$ and $C_2 = B_1^{-1}B_2$. We then remark that if $C_1$
and $C_2$ are diagonalizable matrices, the sum $C_1 \otimes I + I \otimes C_2$ 
can be simultaneously diagonalized by means of the eigendecomposition. More 
generally, the Schur factorization of the matrices $C_1$ and $C_2$ is 
guaranteed to exist, and thus the summation $C_1 \otimes I + I \otimes C_2$ can 
be simultaneously (quasi)-triangularized by of the (real) Schur decomposition.
That is to say, we find orthogonal transformation matrices $Q_1$ and $Q_2$ such 
that
\begin{gather}
  C_1 = Q_1 T_1 Q_1^T, \\
  C_2 = Q_2 T_2 Q_2^T,
\end{gather}
where $T_1$ and $T_2$ are quasi-triangular matrices. Our numerical experiments 
have indicated that the Schur factorizations results in better numerical 
conditioning than the eigendecomposition, and thus we elect to triangularlize 
the matrix rather than diagonalize. Therefore, we can reformulate the linear 
system as
\begin{equation} \label{eq:quad-triang}
\begin{aligned}
  (Q_1 \otimes Q_2) (T_1 \otimes I + I \otimes T_2) (Q_1^T \otimes Q_2^T) x
  &= \left(Q_1 T_1 Q_1^T \otimes I + I \otimes Q_2 T_2 Q_2^T \right) x \\
  &= (C_1 \otimes I + I \otimes C_2) x \\
  &= (A_2^{-1} \otimes B_1^{-1})b.
\end{aligned}
\end{equation}
Since the matrices $Q_1$ and $Q_2$ are orthogonal, the inverse of their 
Kronecker product $Q_1 \otimes Q_2$ is trivially given by $Q_1^T \otimes Q_2^T$.
Thus, solving the system \eqref{eq:quad-triang} is reduced to solving a system 
of the form $T_1 \otimes I + I \otimes T_2$. Well-known solution techniques 
exist for this Sylvester-type system of equations, which can be solved in 
$\mathcal{O}(p^3)$ operations. Thus, once the approximate preconditioner 
$P = A_1 \otimes B_1 + A_2 \otimes B_2$ has been computed, solving linear 
systems of the form $Px = b$ can be performed in linear time per degree of 
freedom.

\subsubsection{Efficient computation of $\widetilde{A}v$ and $\widetilde{A}^Tv$
               in two dimensions} \label{sec:2d-shuffled-matvec}
One of the key operations in efficiently computing the approximate 
Kronecker-product preconditioner is the fast, shuffled matrix-vector product 
operation used in the Lanczos algorithm. Since our algorithm avoids the explicit
construction and evaluation of the entries of the matrix $A$, we present a 
matrix-free algorithm to compute the shuffled product in $\mathcal{O}(p^3)$ 
time. Setting $A = M - \Delta t J$, we apply Algorithm \ref{alg:matvec} to 
compute the shuffled product $\widetilde{A} u$ for a given vector $u$. We first 
write
\begin{equation}
  \widetilde{A} = \widetilde{M} - \Delta t (\widetilde{J_{\rm v}}
    + \widetilde{J_{\rm f}}),
\end{equation}
where $J_{\rm v}$ and $J_{\rm f}$ are the volume and face contributions to the 
Jacobian matrix, respectively. We first demonstrate the computation of the 
product $\widetilde M v$. Recall that the entries of $M$ are given by
\begin{equation} \label{eq:M-entries}
  M_{ij,k\ell} = \int_{K} \tilde\Phi_{ij}(x,y)\tilde\Phi_{k\ell}(x,y)~dx.
\end{equation}
Then, Algorithm \ref{alg:matvec} allows us to write $u = \widetilde{A} v$
\begin{equation}
  u_{:i} = \sum_j (M_{::,ji})^T v_{:j}.
\end{equation}
Writing out the matrix-vector product explicitly, and expanding the integral 
in \eqref{eq:M-entries} as a sum over quadrature nodes, we 
have
\begin{equation}
  u_{ki} = \sum_j \sum_\ell \sum_\alpha \sum_\beta
    \phi_\ell(x_\alpha) \phi_k(x_\beta) \phi_j(x_\alpha) \phi_i(x_\beta)
    |\det(J_T(x_\alpha, x_\beta))|
    w_\alpha w_\beta v_{\ell j}.
\end{equation}
This sum can be factorized as 
\begin{equation}
  u_{ki} = \sum_\beta w_\beta \phi_k(x_\beta) \phi_i(x_\beta)
           \sum_\alpha w_\alpha |\det(J_T(x_\alpha, x_\beta))|
           \sum_j \phi_j(x_\alpha) \sum_\ell \phi_\ell (x_\alpha) v_{\ell j},
\end{equation}
where we notice that each summation in this expression involves no more than 
two free indices, and therefore the expression can be computed in 
$\mathcal{O}(p^3)$ time.

Following the same procedure, and recalling the representation for $J_{\rm v}$ 
and $J_{\rm f}$ given in \eqref{eq:jac-entries}, we can evaluate the shuffled 
product $u = \widetilde{J_{\rm v}} v$ as
\begin{multline}
  u_{ki} = \sum_\beta w_\beta \phi_k(x_\beta) \phi_i(x_\beta)
           \sum_\alpha w_\alpha 
              \frac{\partial F_1}{\partial u}(x_\alpha, x_\beta)
           \sum_j \phi_j(x_\alpha) \sum_\ell \phi'_\ell (x_\alpha) v_{\ell j} \\
    + \sum_\beta w_\beta \phi_k'(x_\beta) \phi_i(x_\beta)
           \sum_\alpha w_\alpha 
              \frac{\partial F_2}{\partial u}(x_\alpha, x_\beta)
           \sum_j \phi_j(x_\alpha) \sum_\ell \phi_\ell (x_\alpha) v_{\ell j}.
\end{multline}
Finally, we consider the face integral terms, and write out the factorized 
form of the shuffled product $u = \widetilde{J_{\rm f}} v$, which takes the 
form
\begin{equation}
  u_{ki} = -\sum_{e \in \partial K} \sum_\alpha w_\alpha
              \phi_i(y_\alpha^e) \phi_k(y_\alpha^e)
              \frac{\partial\widehat{F}}{\partial u^-}(x_\alpha^e, y_\alpha^e)
           \sum_j \phi_j(x_\alpha^e) \sum_\ell \phi_\ell(x_\alpha^e) v_{\ell j}.
\end{equation}
We further remark that many of the terms in this sum can be eliminated by using 
the fact that many of the basis functions are identically zero along a given 
face $e$ of the element $K$.

Computation of the transpose of the shuffled product $\widetilde{A}^T v$ can be 
performed using a very similar matrix-free approach, following the framework 
of Algorithm \ref{alg:matvecT}. These two procedures allow for the computation 
of steps \ref{alg:step-mv} and \ref{alg:step-mvT} in the Lanczos algorithm in 
$\mathcal{O}(p^3)$ time.

\subsection{Three spatial dimensions}
In the case of three spatial dimensions, it would be natural to consider a 
preconditioner matrix $P$ of the form
\begin{equation} \label{eq:3d-precond}
  A \approx P = A_1 \otimes B_1 \otimes C_1 + A_2 \otimes B_2 \otimes C_2 + 
      A_3 \otimes B_3 \otimes C_3.
\end{equation}
Unfortunately, it is not readily apparent how to solve a general system of the 
form \eqref{eq:3d-precond}. Therefore, we instead look for a preconditioner 
that has the simplified form
\begin{equation} \label{eq:3d-simple-precond}
  A \approx P = A_1 \otimes B_1 \otimes C_1 + A_1 \otimes B_2 \otimes C_2,
\end{equation}
where we emphasize that the same matrix $A_1$ appears in both terms on the
right-hand side. This has the advantage that the system $Px = b$ can be
transformed by multiplying on the left by $A_1^{-1} \otimes B_2^{-1} \otimes
C_1^{-1}$ to obtain
\begin{equation}
  (I \otimes B_2^{-1}B_1 \otimes I + I \otimes I \otimes C_1^{-1}C_2)x =
  (A_1^{-1} \otimes B_2^{-1} \otimes C_1^{-1})b.
\end{equation}
Applying the same technique as in the two-dimensional case allows us to 
simultaneously quasi-triangularize both terms on the left-hand side, which 
then results in a system of the form
\begin{equation}
  (I \otimes Q_1 \otimes Q_2)(I \otimes T_1 \otimes I + I \otimes I \otimes T_2)
  (I \otimes Q_1^T \otimes Q_2^T) x
  = (A_1^{-1} \otimes B_2^{-1} \otimes C_1^{-1})b,
\end{equation}
which, as in the case of the two-dimensional system, is a Sylvester-type system 
that can be efficiently solved in $\mathcal{O}(p^3)$ time (constant time 
in $p$ per degree of freedom).

\subsubsection{Forming the three-dimensional preconditioner}
We now address how to generate an effective preconditioner of the form 
\eqref{eq:3d-simple-precond} using the KSVD. First, we recall that the 
element Jacobian will be a $(p+1)^3 \times (p+1)^3$ matrix. We wish to 
approximate this matrix by a Kronecker product $A_1 \otimes D_1$, where 
$A_1 \in \mathbb{R}^{(p+1) \times (p+1)}$ and $D_1 \in \mathbb{R}^{(p+1)^2 
\times (p+1)^2}$. We find such matrices $A_1$ and $D_1$ by finding the largest 
singular value and corresponding singular vectors of the permuted matrix 
$\widetilde{A}$, obtaining 
\begin{equation} \label{eq:3d-step-1}
A \approx A_1 \otimes D_1.
\end{equation}
In order to find the singular values using the Lanczos algorithm, we must
compute the matrix-vector product $\widetilde{A}x$ and $\widetilde{A}^Tx$. In
the following section, we describe how to perform Algorithms \ref{alg:matvec}
and \ref{alg:matvecT} efficiently by taking advantage of the tensor-product
structure of the Jacobian. Once the matrices $A_1$ and $D_1$ have been obtained,
we can then repeat the KSVD process to find the best two-term approximation
\begin{equation} \label{eq:3d-step-2}
  D_1 \approx B_1 \otimes C_1 + B_2 \otimes C_2.
\end{equation}
This too involves the Lanczos algorithm, but since the matrix $D_1$ has 
dimensions $(p+1)^2 \times (p+1)^2$, computing the permuted products 
$\widetilde{D_1}x$ and $\widetilde{D_1}^Tx$ using standard dense linear algebra 
requires $\mathcal{O}(p^4)$ operations, and thus is linear in $p$ per degree 
of freedom. Combining \eqref{eq:3d-step-1} and \eqref{eq:3d-step-2}, we obtain 
an approximation of the form 
\begin{equation}
  A \approx A_1 \otimes B_1 \otimes C_1 + A_1 \otimes B_2 \otimes C_2
\end{equation}
as desired.

\subsubsection{Efficient computation of $\widetilde{A}v$ and 
$\widetilde{A}^Tv$ in three dimensions}
As in Section \ref{sec:2d-shuffled-matvec}, we describe the matrix-free 
procedure for computing the shuffled matrix-vector products 
$\widetilde{A} v$ and $\widetilde{A}^T v$. The general approach to this method 
is the same as in the two-dimensional case, but there are several key 
differences that increase the complexity of this problem. First, we recall that 
since our approximation takes the form $A \approx A_1 \otimes D_1$ where 
$A_1$ is $(p+1) \times (p+1)$ and $D$ is $(p+1)^2 \times (p+1)^2$, the matrix 
$\tilde{A}$ is rectangular, with dimensions $(p+1)^2 \times (p+1)^4$. The 
algorithm we describe has linear complexity per degree of freedom of the vector 
$v$, which results in $\mathcal{O}(p^5)$ operations for the product 
$\widetilde{A} v$, unfortunately not meeting our overall goal of linear time per
degree of freedom in the solution vector.

As before, we decompose the matrix $\widetilde{A} = \widetilde{M} - 
\Delta t (\widetilde{J_{\rm v}} + \widetilde{J_{\rm f}})$. First, we describe 
the method for the mass matrix. Recall that the entries of $M$ are given by
\begin{equation}
  M_{ijk,\ell m n} = \sum_{\alpha,\beta,\gamma}
    w_\alpha w_\beta w_\gamma |\det(J_T(x_\alpha, x_\beta, x_\gamma))| 
    \phi_i(x_\alpha) \phi_j(x_\beta) \phi_k(x_\gamma)
    \phi_\ell(x_\alpha) \phi_m(x_\beta) \phi_n(x_\gamma)~dx.
\end{equation}
Then, following Algorithm \ref{alg:matvec}, we write $u = \widetilde{M}v$, 
where $u$ is a vector of length $(p+1)^2$ and $v$ is a vector of length 
$(p+1)^4$,
\begin{equation}
  u_{:i} = \sum_j \sum_k (M_{:::,jki})^T v_{::jk}.
\end{equation}
Following the same sum factorization procedure as in the two-dimensional case,
we can write
\begin{multline}
  u_{\ell i} = \sum_\gamma \phi_i(x_\gamma) \phi_\ell(x_\gamma)
               \sum_\beta 
               \sum_\alpha w_\alpha w_\beta w_\gamma
                  | \det(J_T(x_\alpha, x_\beta, x_\gamma)) | \\
               \sum_k \phi_k(x_\beta)
               \sum_j \phi_j(x_\alpha)
               \sum_n \phi_n(x_\beta)
               \sum_m \phi_m(x_\alpha) v_{mnjk}.
\end{multline}
The $\mathcal{O}(p^5)$ complexity is clear from this form, as, for example, 
the right-most summation has four free indices. The shuffled product with the 
Jacobian of the volume integral takes a similar form, where 
$u = \widetilde{J_{\rm v}} v$ can be written as
\begin{multline}
  u_{\ell i} = \sum_\gamma \phi_i(x_\gamma) \phi_\ell(x_\gamma) 
               \sum_\beta \sum_\alpha  w_\alpha w_\beta w_\gamma \\
           \bigg( \frac{\partial F_1}{\partial u}(x_\alpha, x_\beta, x_\gamma)
            \sum_k \phi_k(x_\beta) \sum_j \phi_j(x_\alpha)
            \sum_n \phi_n(x_\beta) \sum_m \phi_m'(x_\alpha) v_{mnjk} \\
              + \frac{\partial F_2}{\partial u}(x_\alpha, x_\beta, x_\gamma)
            \sum_k \phi_k(x_\beta) \sum_j \phi_j(x_\alpha)
            \sum_n \phi_n'(x_\beta) \sum_m \phi_m(x_\alpha) v_{mnjk} \bigg) \\
            + \sum_\gamma \phi_i(x_\gamma) \phi_\ell'(x_\gamma) 
               \sum_\beta \sum_\alpha  w_\alpha w_\beta w_\gamma
               \frac{\partial F_3}{\partial u}(x_\alpha, x_\beta, x_\gamma)\\
               \sum_k \phi_k(x_\beta) \sum_j \phi_j(x_\alpha)
               \sum_n \phi_n(x_\beta) \sum_m \phi_m(x_\alpha) v_{mnjk}.
\end{multline}
Finally, we write the product corresponding to the face integral Jacobian, 
$u = \widetilde{J_{\rm f}}v$, as
\begin{multline}
  u_{\ell i} = \sum_{e \in \partial K} \sum_\beta \sum_\alpha w_\alpha w_\beta
                  \phi_i(z_{\alpha\beta}^e)
                  \phi_\ell(z_{\alpha\beta}^e) 
                  \frac{\partial \widehat{F}}{\partial u^-}
                  (x_{\alpha\beta}^e,y_{\alpha\beta}^e,z_{\alpha\beta}^e) \\
               \sum_k \phi_k(y_{\alpha\beta}^e)
               \sum_j \phi_j(x_{\alpha\beta}^e)
               \sum_n \phi_n(y_{\alpha\beta}^e)
               \sum_m \phi_m(x_{\alpha\beta}^e) v_{mnjk},
\end{multline}
where $(x_{\alpha\beta}^e, y_{\alpha\beta}^e, z_{\alpha\beta}^e)$ represents the
coordinates of the quadrature nodes on the face $e$ of element $K$ indexed by
$(\alpha, \beta)$. We remark that for each face $e$, two of the indices in the
above expression can be eliminated. This simplification depends on the 
orientation of the face, and therefore we leave the full expression for the sake
of generality.

\subsection{Algorithm overview}
Here we describe the overall algorithms used to form and apply the tensor
product preconditioner. We present the algorithm for both the cases of two and
three spatial dimensions. Forming the preconditioner requires the Lanczos SVD,
given by Algorithm \ref{alg:lanczos}, and the two permuted matrix-vector
multiplication kernels, given by Algorithms \ref{alg:matvec} and
\ref{alg:matvecT}, and described in the preceding section. Computational
complexities are indicated for each step of the algorithm. We note that in the
2D case, we obtain an overall complexity of $\mathcal{O}(p^{d+1})$. In the 3D
case, all the operations have complexity at most $\mathcal{O}(p^{d+1})$, except
the first Lanczos SVD, which requires $\mathcal{O}(p^5)$ operations. 

\begin{algorithm}[H]
\caption{Form 2D preconditioner $A \approx P = A_1 \otimes B_2 
          + A_2 \otimes B_2$}
\label{alg:2d-form}
\begin{algorithmic}[1]
  \State Compute $A \approx A_1 \otimes B_1 + A_2 \otimes B_2$ using Lanczos
         iteration and matrix-free products
         $\widetilde{A}x$ and $\widetilde{A}^T x$
         \Comment{Complexity: $\mathcal{O}(p^3)$}
  \State Precompute LU factorizations of $A_2$ and $B_1$
         \Comment{Complexity: $\mathcal{O}(p^3)$}
  \State Precompute Schur factorizations $Q_1 T_1 Q_1^T, Q_2 T_2 Q_2^T$
         of $A_2^{-1} A_1$ and $B_1^{-1} b_2$, respectively
         \Comment{Complexity: $\mathcal{O}(p^3)$}
\end{algorithmic}
\end{algorithm}
\begin{algorithm}[H]
\caption{Apply 2D preconditioner to solve $Px = b$}
\label{alg:2d-apply}
\begin{algorithmic}[1]
  \State $\tilde{b} \gets A_2^{-1} \otimes B_1^{-1} b$
         \Comment{Complexity: $\mathcal{O}(p^3)$}
  \State Solve the Sylvester system
         $(T_1 \otimes I + I \otimes T_2) \tilde{x}
          = (Q_1^T \otimes Q_2^T) \tilde{b}$
         \Comment{Complexity: $\mathcal{O}(p^3)$}
  \State $x \gets (Q_1 \otimes Q_2) \tilde{x}$
         \Comment{Complexity: $\mathcal{O}(p^3)$}
\end{algorithmic}
\end{algorithm}
\begin{algorithm}[H]
\caption{Form 3D preconditioner $A \approx P = A_1 \otimes B_1 \otimes C_1 
         + A_1 \otimes B_2 \otimes C_2$}
\label{alg:3d-form}
\begin{algorithmic}[1]
  \State Compute $A \approx A_1 \otimes D_1$ using Lanczos iteration and 
         matrix-free products $\widetilde{A}x$ and $\widetilde{A}^T x$
         \Comment{Complexity: $\mathcal{O}(p^5)$}
  \State Compute $D_1 \approx B_1 \otimes C_1 + B_2 \otimes C_2$ using 
         Lanczos iteration and dense permuted products
         \Comment{Complexity: $\mathcal{O}(p^4)$}
  \State Precompute LU factorizations of $A_1$, $B_2$, and $C_1$
         \Comment{Complexity: $\mathcal{O}(p^3)$}
  \State Precompute Schur factorizations $Q_1 T_1 Q_1^T, Q_2 T_2 Q_2^T$
         of $B_2^{-1} B_1$ and $C_1^{-1} C_2$, respectively
         \Comment{Complexity: $\mathcal{O}(p^3)$}
\end{algorithmic}
\end{algorithm}
\begin{algorithm}[H]
\caption{Apply 3D preconditioner to solve $Px = b$}
\label{alg:3d-apply}
\begin{algorithmic}[1]
  \State $\tilde{b} \gets A_1^{-1} \otimes B_2^{-1} \otimes C_1^{-1} b$
         \Comment{Complexity: $\mathcal{O}(p^4)$}
  \State Solve the Sylvester system
    $(I \otimes T_1 \otimes I + I \otimes I \otimes T_2) \tilde{x}
         = (I \otimes Q_1^T \otimes Q_2^T) \tilde{b}$
         \Comment{Complexity: $\mathcal{O}(p^4)$}
  \State $x \gets (I \otimes Q_1 \otimes Q_2) \tilde{x}$
         \Comment{Complexity: $\mathcal{O}(p^4)$}
\end{algorithmic}
\end{algorithm}

%%%%%%%%%%%%%%%%%%%%%%%%%%%%%%%%%%%%%%%%%%%%%%%%%%%%%%%%%%%%%%%%%%%%%%%%%%%%%%%%
% Numerical results
%%%%%%%%%%%%%%%%%%%%%%%%%%%%%%%%%%%%%%%%%%%%%%%%%%%%%%%%%%%%%%%%%%%%%%%%%%%%%%%%

\section{Numerical results}\label{sec:results}

{
In the following sections, we present numerical results that demonstrate several
important features of the preconditioner and its performance when applied to 
a variety of equations and test cases. We consider both two-dimensional and 
three-dimensional problems, and solve the scalar advection equation, the Euler 
equations, and the Navier-Stokes equations. The nonlinear systems of equations 
resulting from implicit time integration are solved using Newton's method, with 
a relative tolerance of $10^{-8}$. Within each Newton iteration, the linear 
system is solved using preconditioned, restarted GMRES, with a relative 
tolerance of $10^{-5}$. The parameters of Newton tolerance, GMRES tolerance, and
restart iterations are chosen according to performance study found in 
\cite{Zahr2013}. Although these parameters can have an effect on overall 
solution time, the relationships are often neither simple nor well-understood, 
and these issues are not considered in depth in this work.
}

\subsection{2D linear advection equation}
The simplest example we consider is that of the two-dimensional scalar advection
equation, given by
\begin{equation} \label{eq:adv}
  u_t + \nabla \cdot \left( \alpha, \beta \right) u = 0,
\end{equation}
where $(\alpha, \beta)$ is a space-dependent velocity field. Because of the 
particularly simple structure of this equation, it is possible to see how the 
approximate Kronecker preconditioner, given by $A_1 \otimes B_1 + 
A_2 \otimes B_2$ relates to the true diagonal blocks $A$ of the matrix 
$M - \Delta t J$. In this case, the properties of the velocity field 
$(\alpha, \beta)$ can determine how well the discontinuous Galerkin 
discretization can be approximated by a tensor-product structure.

Neglecting for now the face integral terms, the diagonal blocks of 
$M - \Delta t J$ can be written as
\begin{equation} \label{eq:adv-kron}
  \left( G^T W \otimes G^T W \right) J_T \left( G \otimes G \right)
   - \Delta t \left(
   \left( G^T W \otimes D^T W \right) F_1
   +
   \left( D^T W \otimes G^T W \right) F_2 
   \right) \left( G \otimes G \right),
\end{equation}
where $J_T, F_1,$ and $F_2$ are $\mu^2 \times \mu^2$ diagonal matrices. If the 
matrices $J_T, F_1,$ and $F_2$ additionally posses a Kronecker-product 
structure, then it is possible to rewrite \eqref{eq:adv-kron} exactly in the 
form $A_1 \otimes B_1 + A_2 \otimes B_2$. The Kronecker structure of $J_T$ is 
determined by the geometry of the mesh, and the structure of $F_1$ and $F_2$ 
is determined by the form of the velocity field.

For example, we first suppose that the mesh is a cartesian grid with grid size
$h$, and thus the Jacobian determinant of the transformation map is equal to
$h^2$. Hence, $J_T$ is equal to $h^2$ times the identity matrix. If we further
suppose that the velocity field is separable, in the sense that, each component
depends only on the corresponding spatial variable, \textit{i.e.} $\alpha(x,y) =
\alpha(x)$, $\beta(x,y) = \beta(y)$, then the flux derivatives can be written as
$F_1 = I \otimes F_{1_x}$, and $F_2 = F_{2_y} \otimes I$. Therefore, we can
rewrite \eqref{eq:adv-kron} in the form
\begin{equation}  \label{eq:adv-kron-cart}
  \left( h^2 G^TWG - \Delta t D^T W F_{2_y} G \right) \otimes G^T W G
  - \Delta t \left( G^T W G \otimes D^T W F_{1_x} G \right),
\end{equation}
and we see that the diagonal blocks of $M - \Delta t J$ are exactly 
representable by our Kronecker-product approximation.

If, on the other hand, we allow straight-sided, non-cartesian meshes, then the 
transformation mapping is a bilinear function, and its Jacobian determinant is a
linear function in the variables $x$ and $y$. Thus, $J_T = J_{T_y} \otimes I 
+ I \otimes J_{T_x}$. If the velocity field is constant in space, then we obtain
the following representation of the diagonal blocks
\begin{equation} \label{eq:adv-kron-const}
  G^T W G \otimes \left( G^T W J_{T_x} G - \Delta t \alpha D^T W G \right) +
  \left( G^T W J_{T_y} G - \Delta t \beta D^T W G \right) \otimes G^T W G,
\end{equation}
and we see that our Kronecker-product approximation is again exact.

If, in contrast to the previous two cases, the transformation mapping is given 
by a higher degree polynomial, or if the velocity field is not separable, then 
the approximate preconditioner will not yield the exact diagonal blocks. 
However, if the deformation of the mesh is not too large, if the velocity 
field is well approximated by one that is separable, or if the time step 
$\Delta t$ is relatively small, then we expect the Kronecker product 
preconditioner to compare favorably with the exact block Jacobi preconditioner. 

\begin{figure}[t!]
    \centering
    \begin{subfigure}[t]{0.5\textwidth}
        \centering
        \includegraphics[width=2in]{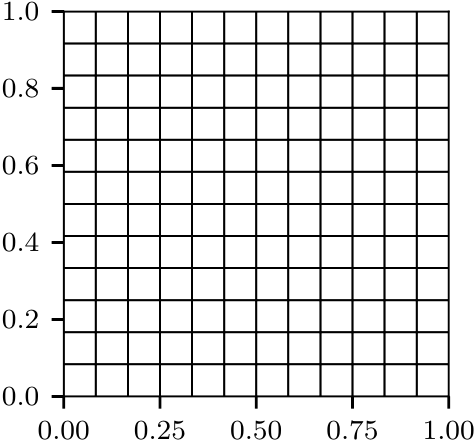}
        \caption{Cartesian mesh for advection equation}
    \end{subfigure}%
    ~ 
    \begin{subfigure}[t]{0.5\textwidth}
        \centering
        \includegraphics[width=2in]{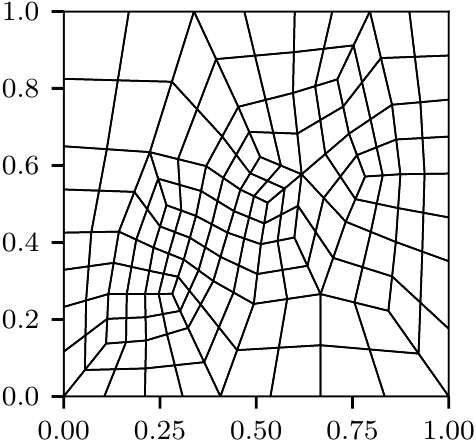}
        \caption{Unstructured mesh for advection equation}
        \label{fig:adv-unstructured-mesh}
    \end{subfigure}
    \caption{Meshes used for the advection equation}
    \label{fig:adv-meshes}
\end{figure}
\begin{figure}[t!]
    \centering
    \begin{subfigure}[t]{0.32\textwidth}
        \centering
        \includegraphics[width=2in]{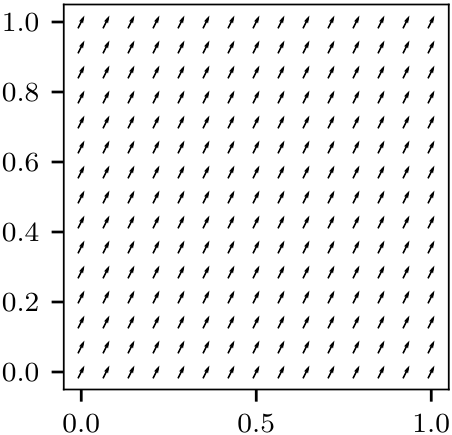}
        \caption{$\alpha = 1, \beta = 2$} \label{fig:vf1}
    \end{subfigure}%
    ~ 
    \begin{subfigure}[t]{0.32\textwidth}
        \centering
        \includegraphics[width=2in]{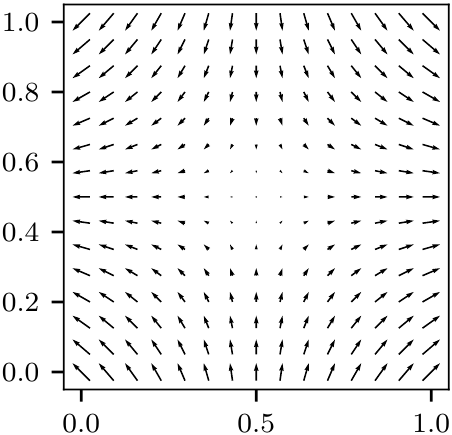}
        \caption{$\alpha = x - \frac{1}{2}, \beta = \frac{1}{2} - y$}
        \label{fig:vf2}
    \end{subfigure}
    ~ 
    \begin{subfigure}[t]{0.32\textwidth}
        \centering
        \includegraphics[width=2in]{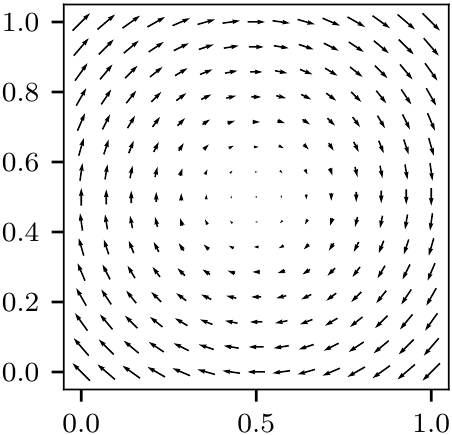}
        \caption{$\alpha = y - \frac{1}{2}, \beta = \frac{1}{2} - x$}
        \label{fig:vf3}
    \end{subfigure}
    \caption{Velocity fields used for the advection equation}
    \label{fig:adv-vel-field}
\end{figure}

We remark that this numerical experiment is designed to highlight two main 
features of the Kronecker-product preconditioner. The first is that if any
exact representation of the diagonal blocks in the form $A_1 \otimes B_1 + A_2
\otimes B_2$ exists, such as those given by equations \eqref{eq:adv-kron-cart}
and \eqref{eq:adv-kron-const}, then the KSVD algorithm provides an
\textit{automatic}, and purely \textit{algebraic} method to identity this
decomposition. No special structure of the flux functions is required to be
known a priori in order for the KSVD to exactly reproduce this tensor-product
structure. Secondly, even in a case where it is impossible to write such an
expression exactly, the KSVD method will identify the best possible
approximation of this form. Thus, in cases where the velocity field is close to
constant (\textit{e.g.} when the mesh size is very small), or where the mesh
deformation is small, we expect this approximation to be very accurate.

In order to compare the performance of these two preconditioners, we solve 
equation \eqref{eq:adv} on both regular and irregular meshes, with constant, 
separable, and non-separable velocity fields. The meshes are shown in Figure 
\ref{fig:adv-meshes}, and the velocity fields in Figure \ref{fig:adv-vel-field}.
We choose a representative time step of $\Delta t = 0.5$, and consider 
polynomial degrees $p = 1, 2, \ldots, 10$. In the case of the regular cartesian 
mesh, we expect identical performance for the exact block Jacobi and approximate
Kronecker-product preconditioners for the velocity fields shown in Figures 
\ref{fig:vf1} and \ref{fig:vf2}, since the diagonal blocks can be reproduced 
exactly. In the case of the unstructured mesh, we expect to see identical 
performance for the constant velocity field in Figure \ref{fig:vf1}. Indeed, 
the numerical results corroborate our expectations, and the number of iterations
is identical between the two preconditioners in those test cases. Additionally,
even in cases where the Kronecker-product approximation cannot reproduce the 
exact blocks, such for for velocity field \ref{fig:vf3} or \ref{fig:vf2} on the 
unstructured mesh, the performance is, in most cases, extremely similar to that 
of exact block Jacobi. The number of GMRES iterations required to converge 
with each preconditioner is shown in Table \ref{tab:adv-results}. For very large
values of polynomial degree $p$ on the unstructured mesh, with non-separable 
velocity field, we begin to see a degradation in the performance of the 
Kronecker-product preconditioner.

\begin{table}[H]
\caption{Number of GMRES iterations for Jacobi (\textit{J}) and KSVD 
         (\textit{K}) preconditioners, advection equation on cartesian and 
         unstructured grids, velocity fields (a), (b), and (c) from Figure 
         \ref{fig:adv-vel-field}.}
\label{tab:adv-results}
\centering
\begin{subtable}{0.4\linewidth}
\caption{Cartesian grid}
\begin{tabular}{l|ll|ll|ll}
\toprule
& \multicolumn{2}{c|}{\textbf{(a)}}
& \multicolumn{2}{c|}{\textbf{(b)}}
& \multicolumn{2}{c}{\textbf{(c)}}\\
$p$&\textit{J}&\textit{K}&\textit{J}&\textit{K}&\textit{J}&\textit{K} \\
\midrule
1 & 12 & 12 & 5 & 5 & 29 & 29 \\
2 & 14 & 14 & 7 & 7 & 29 & 29 \\
3 & 13 & 13 & 7 & 7 & 29 & 29 \\
4 & 14 & 14 & 7 & 7 & 29 & 29 \\
5 & 13 & 13 & 7 & 7 & 29 & 30 \\
6 & 17 & 17 & 7 & 7 & 29 & 31 \\
7 & 13 & 13 & 7 & 7 & 30 & 28 \\
8 & 14 & 14 & 7 & 7 & 30 & 29 \\
9 & 12 & 12 & 7 & 7 & 27 & 30 \\
10 & 14 & 14 & 7 & 7 & 27 & 30 \\
\bottomrule
\end{tabular}
\end{subtable}
\begin{subtable}{0.4\linewidth}
\caption{Unstructured mesh}
\label{tab:adv-unstructured}
\begin{tabular}{l|ll|ll|ll}
\toprule
& \multicolumn{2}{c|}{\textbf{(a)}}
& \multicolumn{2}{c|}{\textbf{(b)}}
& \multicolumn{2}{c}{\textbf{(c)}}\\
$p$&\textit{J}&\textit{K}&\textit{J}&\textit{K}&\textit{J}&\textit{K} \\
\midrule
1 & 14 & 14 & 10 & 11 & 29 & 29 \\
2 & 15 & 15 & 11 & 10 & 29 & 29 \\
3 & 14 & 14 & 11 & 12 & 28 & 28 \\
4 & 15 & 15 & 9 & 12 & 28 & 31 \\
5 & 14 & 14 & 10 & 12 & 28 & 34 \\
6 & 14 & 14 & 10 & 12 & 28 & 39 \\
7 & 13 & 13 & 10 & 12 & 28 & 46 \\
8 & 13 & 13 & 12 & 13 & 28 & 53 \\
9 & 13 & 13 & 12 & 14 & 28 & 62 \\
10 & 13 & 13 & 12 & 15 & 28 & 69 \\
\bottomrule
\end{tabular}
\end{subtable}
\end{table}

{
\subsection{Anisotropic grids}

One main motivation for the use of implicit time integration methods is the
presence of stretched or highly anisotropic elements, for instance in the
vicinity of a shock, or at a boundary layer \cite{yano2012optimization}. In order to investigate
the performance of the Kronecker-product preconditioner for this important class
of problems, we consider the scalar advection equation on two anisotropic grids,
shown in Figure \ref{fig:aniso-meshes}. The first mesh consists entirely of
rectangular elements, refined around the center line $x = 1/2$, such that the
thinnest elements have an aspect ratio of about 77. The second mesh is similar,
with the main difference being that the quadrilateral elements no longer posses
$90^\circ$ angles. In accordance with the analysis from the preceding section,
we can expect the Kronecker-product preconditioner to exactly reproduce the
diagonal blocks in the rectangular case for separable velocity fields. However,
in the case of the skewed quadrilaterals, the Kronecker preconditioner is only
exact for constant velocity fields, and provides an approximation to the
diagonal blocks of the Jacobian in other cases. In the interest of generality,
we consider the non-separable velocity field (c) shown in Figure
\ref{fig:adv-vel-field}, for which the Kronecker-product preconditioner is 
approximate for both the rectangular and skewed meshes.

\begin{figure}[t!]
    \centering
    \begin{subfigure}[t]{0.5\textwidth}
        \centering
        \includegraphics[width=2in]{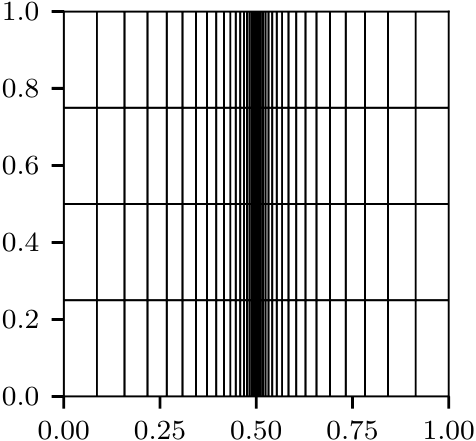}
        %\caption{Cartesian mesh for advection equation}
    \end{subfigure}%
    ~ 
    \begin{subfigure}[t]{0.5\textwidth}
        \centering
        \includegraphics[width=2in]{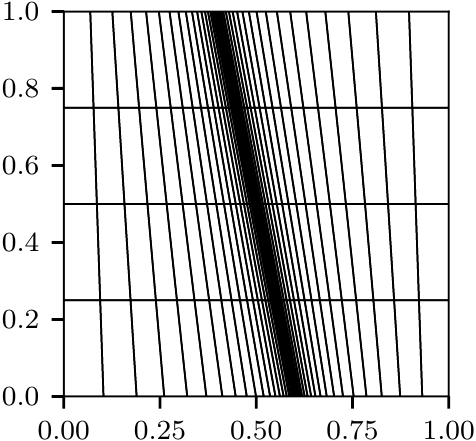}
        %\caption{Unstructured mesh for advection equation}
    \end{subfigure}
    \caption{ Meshes used for anisotropic test case}
    \label{fig:aniso-meshes}
\end{figure}

We use this test case to compare the runtime performance of the
Kronecker-product preconditioner both with explicit time integration methods,
and with the exact block Jacobi preconditioner. To this end, we choose an
implicit time step of $\Delta t = 5 \times 10^{-2}$. We then compute one time
step using a third-order $L$-stable DIRK method \cite{Alexander:1977dk}.
Additionally, we integrate until $t = 5 \times 10^{-2}$ using the standard
fourth-order explicit Runge-Kutta method, with the largest possible stable
explicit time step. The explicit time step restriction becomes more severe as
the polynomial degree $p$ increases \cite{Krivodonova:2013}, resulting in a
large increase in the number of time steps required.

We choose polynomial degrees $p = 1, 2, \ldots, 30$, and measure the runtime
required to integrate until $t = 5 \times 10^{-2}$. Due to the excessive
runtimes, we use only $p= 1, 2, \ldots, 15$ for the explicit method. We display
the runtimes for both rectangular and general quadrilateral meshes in Figure
\ref{fig:aniso-runtimes}. For $p > 1$, the explicit RK4 method is not 
competitive for this problem. For both meshes, the KSVD preconditioner results 
in faster runtimes than the exact block Jacobi preconditioner starting at about
$p = 4$ or $p = 5$. In the rectangular case, we see a noticeable asymptotic 
improvement in the runtime in this case. For $p=30$, the Kronecker-product 
preconditioner results in runtimes close to 20 times faster than block Jacobi.
In the case of the skewed quadrilateral mesh, we observe an increase in the 
number of GMRES iterations required per time step, similar to what was observed 
in column (c) of Table \ref{tab:adv-unstructured}. Despite this increase in 
iteration count, the Kronecker-product preconditioner still resulted in runtimes
about three times shorter than the exact block Jacobi.

\begin{figure}[b!]
    \centering
    \begin{subfigure}[t]{0.5\textwidth}
        \centering
        \includegraphics[width=3in]{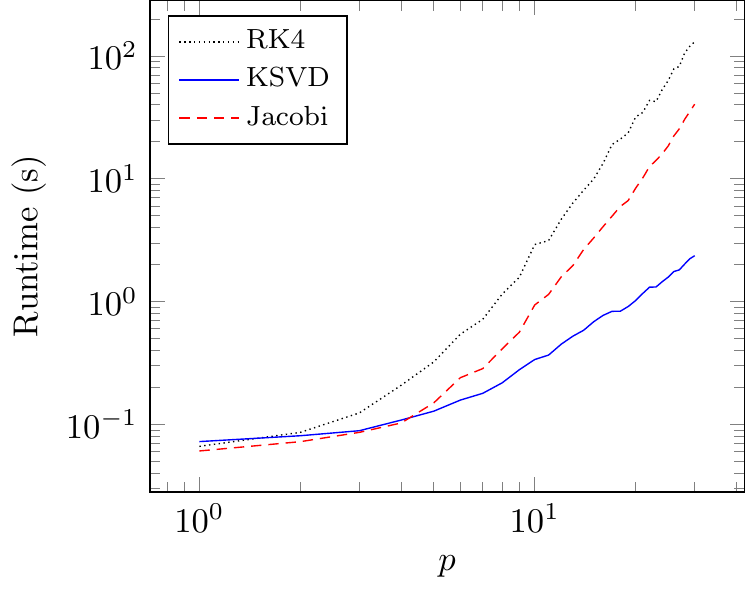}
        \caption{Rectangular anisotropic mesh}
    \end{subfigure}%
    ~ 
    \begin{subfigure}[t]{0.5\textwidth}
        \centering
        \includegraphics[width=3in]{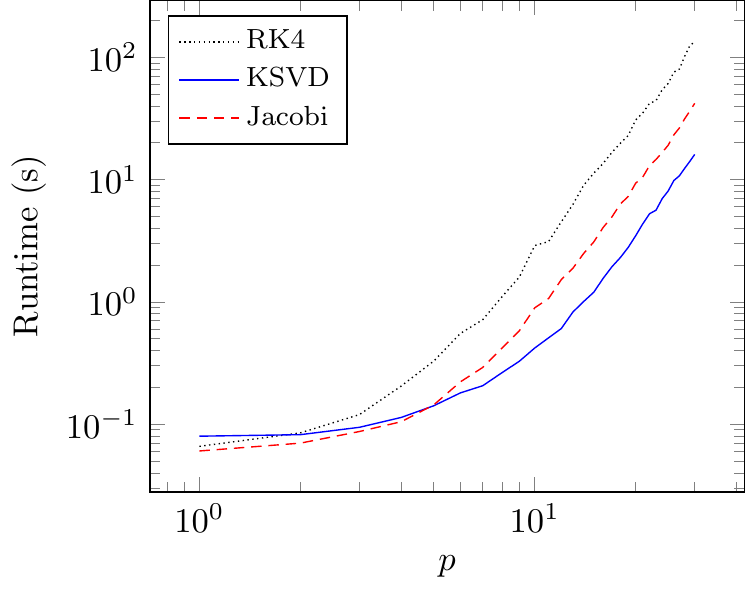}
        \caption{Skewed anisotropic mesh}
    \end{subfigure}
    \caption{Wall-clock times for scalar advection on anisotropic 
             meshes}
    \label{fig:aniso-runtimes}
\end{figure}

Additionally, we measure the average wall-clock time required to both form and
apply the Kronecker and block Jacobi preconditioners, for all polynomial degrees
considered. We see that the cost of forming the Jacobi preconditioner quickly 
dominates the runtime. For large $p$, we begin to see the asymptotic 
$\mathcal{O}(p^6)$ complexity for this operation. Applying the Jacobi 
preconditioner requires $\mathcal{O}(p^4)$ operations, while both forming and 
applying the Kronecker preconditioner require $\mathcal{O}(p^3)$ operations. 
These computational complexities are evident from the measured wall-clock times,
shown in Figure \ref{fig:aniso-form-apply}.

\begin{figure}[t!]
    \centering
    \includegraphics[width=3in]{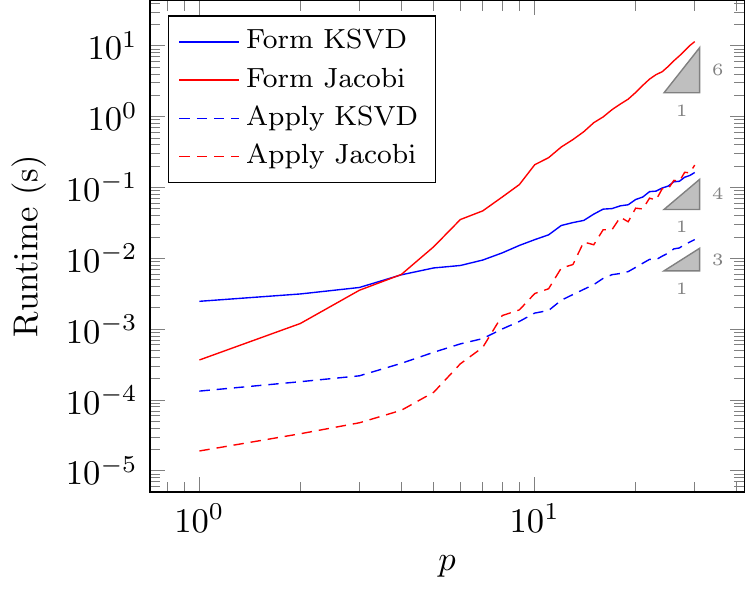}
    \caption{Wall-clock time required to form (solid lines) and 
             apply (dashed lines) the exact block Jacobi and approximate 
             Kronecker-product preconditioners.}
    \label{fig:aniso-form-apply}
\end{figure}
}

\subsection{2D Euler vortex}
In this example, we consider the compressible Euler equations of gas dynamics 
in two dimensions, given in conservative form by 
\begin{equation}
\partial_t \bm{u} + \nabla \cdot F(\bm{u}) = 0,
\end{equation}
where
\begin{equation} \label{eq:euler-flux}
  \def\arraystretch{1}
  \bm{u} = \left(\begin{array}{c}
    \rho \\ \rho u \\ \rho v \\ \rho E
  \end{array}\right), \qquad
  F_1(\bm{u}) = \left(\begin{array}{c}
    \rho u \\ \rho u^2 + p \\ \rho uv \\ \rho Hu
  \end{array}\right), \qquad
  F_2(\bm{u}) = \left(\begin{array}{c}
    \rho v \\ \rho uv \\ \rho v^2 + p \\ \rho Hv
  \end{array}\right),
\end{equation}
where $\rho$ is the density, $\bm{v} = (u, v)$ is the fluid velocity, $p$ is the
pressure, and $E$ is the specific energy. The total enthalpy $H$ is given by 
\begin{equation}
    H = E + \frac{p}{\rho},
\end{equation}
and the pressure is determined by the equation of state
\begin{equation}
    p = (\gamma - 1)\rho \left(E - \frac{1}{2}(u^2 + v^2)\right),
\end{equation}
where $\gamma = c_p/c_v$ is the ratio of specific heat capacities at constant
pressure and constant volume, taken to be 1.4.

We consider the model problem of an unsteady compressible vortex in a
rectangular domain \cite{Wang2013}. The domain is taken to be a $20\times15$
rectangle and the vortex is initially centered at $(x_0, y_0) = (5, 5)$. The
vortex is moving with the free-stream at an angle of $\theta$. The exact
solution is given by
\begin{gather}
    u = u_\infty \left( \cos(\theta) - 
        \frac{\epsilon ((y-y_0) - \overline{v} t)}{2\pi r_c} 
        \exp\left( \frac{f(x,y,t)}{2} \right) \right),\\
    u = u_\infty \left( \sin(\theta) - \frac{\epsilon ((x-x_0) - 
        \overline{u} t)}{2\pi r_c} 
        \exp\left( \frac{f(x,y,t)}{2} \right) \right),\\
    \rho = \rho_\infty \left( 1 - 
        \frac{\epsilon^2 (\gamma - 1)M^2_\infty}{8\pi^2} \exp(f(x,y,t))
        \right)^{\frac{1}{\gamma-1}}, \\
    p = p_\infty \left( 1 - 
        \frac{\epsilon^2 (\gamma - 1)M^2_\infty}{8\pi^2} \exp(f(x,y,t))
        \right)^{\frac{\gamma}{\gamma-1}},
\end{gather}
where $f(x,y,t) = (1 - ((x-x_0) - \overline{u}t)^2 - ((y-y_0) -
\overline{v}t)^2)/r_c^2$, $M_\infty$ is the Mach number, $u_\infty,
\rho_\infty,$ and $p_\infty$ are the free-stream velocity, density, and
pressure, respectively. The free-stream velocity is given by $(\overline{u},
\overline{v}) = u_\infty (\cos(\theta), \sin(\theta))$. The strength of the
vortex is given by $\epsilon$, and its size is $r_c$. We choose the parameters
to be $M_\infty = 0.5$, $u_\infty = 1$, $\theta = \arctan(1/2)$, $\epsilon =
0.3$, and $r_c = 1.5$.

\begin{figure}[t!]
  \centering
  \includegraphics[width=3.0in]{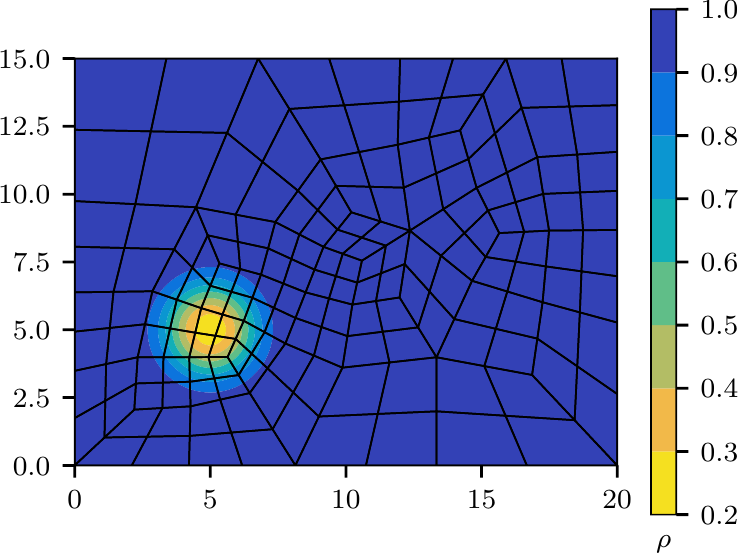}
  \caption{Initial conditions (density) for Euler vortex on unstructured mesh.}
  \label{fig:ev-ic}
\end{figure}

As in the case of the linear advection equation, we consider both a regular 
$n_x \times n_y$ cartesian grid, and an unstructured mesh. The unstructured 
mesh is obtained by scaling the mesh in Figure \ref{fig:adv-unstructured-mesh} 
by 20 in the $x$-direction and 15 in the $y$-direction. Density contours of the 
initial conditions are shown in Figure \ref{fig:ev-ic}. As opposed to the 
scalar advection equation, the solution to the Euler equations consists of 
multiple components. Thus, the blocks of the Jacobian matrix can be considered 
to be of size $n_c (p+1)^2 \times n_c (p+1)^2$, where $n_c$ is the number of 
solution components (in the case of the 2D Euler equations, $n_c = 4$). The 
exact block Jacobi preconditioner computes the inverses of these large blocks. 
The approximate Kronecker-product preconditioner find optimal approximation of 
the form $A_1 \otimes B_1 + A_2 \otimes B_2$, where $A_1$ and $A_2$ are 
$n_c(p+1) \times n_c(p+1)$ matrices, and $B_1$ and $B_2$ are $(p+1)\times(p+1)$
matrices.

We choose two representative time steps of $\Delta t = 0.1$ and $\Delta t=0.01$,
and compute the average number of GMRES iterations per linear solve required to
perform one backward Euler time step. We choose the
polynomial degree $p = 3,4,\ldots,15$, and consider cartesian and unstructured
meshes, both with 160 quadrilateral elements. We present the results in Table
\ref{tab:ev-results}. Very similar results are observed for the structured and 
unstructured results. We note that for the smaller time step, the approximate
Kronecker-product preconditioner requires a very similar number of iterations
when compared with the exact block Jacobi preconditioner, even for high
polynomial degree $p$. For the larger time step, the number of iterations
required for the KSVD preconditioner increases with $p$ at a faster rate when
compared with the block Jacobi preconditioner, suggesting that the 
Kronecker-product preconditioner is most effective for moderate time steps
$\Delta t$.

\begin{table}[b!]
%\small
\caption{Number of GMRES iterations for Jacobi (\textit{J}) and KSVD 
         (\textit{K}) preconditioners, Euler equations on cartesian and 
         unstructured grids, with $\Delta t = 0.1$ and $\Delta t = 0.01$.}
\label{tab:ev-results}
\centering
\begin{subtable}{0.4\linewidth}
\centering
\caption{Cartesian grid}
\begin{tabular}{l|ll|ll}
\toprule
& \multicolumn{2}{c|}{$\Delta t = 0.01$}
& \multicolumn{2}{c}{$\Delta t = 0.1$}\\
$p$&\textit{J}&\textit{K}&\textit{J}&\textit{K} \\
\midrule
%1 & 5 & 5 & 8 & 8 \\
%2 & 5 & 5 & 9 & 13 \\
3 & 5 & 6 & 11 & 18 \\
4 & 6 & 7 & 12 & 23 \\
5 & 6 & 8 & 13 & 30 \\
6 & 7 & 9 & 15 & 38 \\
7 & 7 & 10 & 17 & 47 \\
8 & 8 & 11 & 18 & 59 \\
9 & 8 & 13 & 20 & 71 \\
10 & 9 & 15 & 21 & 88 \\
11 & 9 & 17 & 23 & 103 \\
12 & 10 & 19 & 25 & 121 \\
13 & 11 & 20 & 24 & 123 \\
14 & 11 & 23 & 25 & 157 \\
15 & 12 & 25 & 26 & 196 \\
\bottomrule
\end{tabular}
\end{subtable}
\begin{subtable}{0.4\linewidth}
\centering
\caption{Unstructured mesh}
\begin{tabular}{l|ll|ll}
\toprule
& \multicolumn{2}{c|}{$\Delta t = 0.01$}
& \multicolumn{2}{c}{$\Delta t = 0.1$}\\
$p$&\textit{J}&\textit{K}&\textit{J}&\textit{K} \\
\midrule
%1 & 5 & 5 & 9 & 10 \\
%2 & 5 & 6 & 10 & 14 \\
3 & 6 & 7 & 12 & 19 \\
4 & 6 & 7 & 14 & 26 \\
5 & 7 & 9 & 16 & 32 \\
6 & 7 & 10 & 17 & 42 \\
7 & 8 & 11 & 18 & 51 \\
8 & 8 & 12 & 20 & 64 \\
9 & 9 & 13 & 21 & 74 \\
10 & 9 & 15 & 23 & 90 \\
11 & 9 & 17 & 24 & 110 \\
12 & 10 & 19 & 25 & 125 \\
13 & 10 & 21 & 25 & 142 \\
14 & 11 & 24 & 26 & 164 \\
15 & 11 & 26 & 27 & 245 \\
\bottomrule
\end{tabular}
\end{subtable}
\end{table}

\subsubsection{Performance comparison}
\begin{figure}[h!]
    \centering
    \begin{subfigure}[t]{0.5\textwidth}
        \centering
        \includegraphics[width=2.9in]{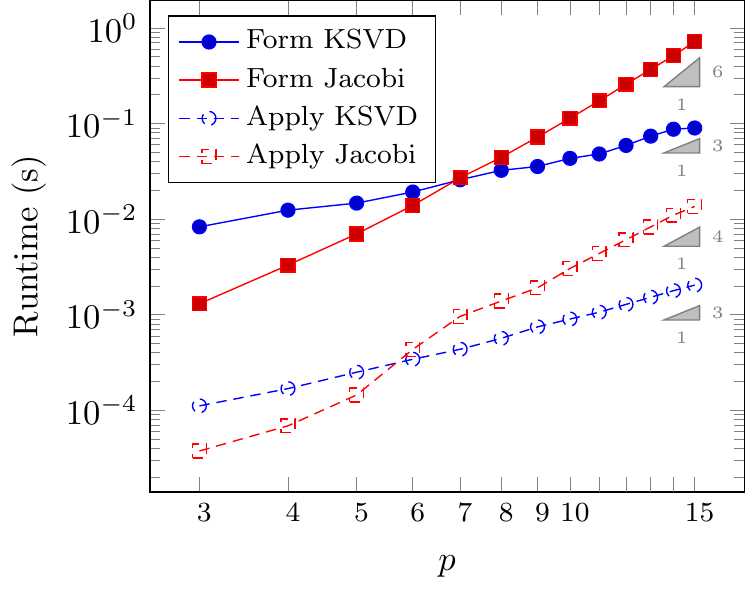}
        \caption{Wall-clock time required to form and apply the KSVD and Jacobi 
                 preconditioners}
        \label{fig:ev-timings}
    \end{subfigure}%
    ~ \hspace{0.1in}
    \begin{subfigure}[t]{0.5\textwidth}
        \centering
        \includegraphics[width=2.9in]{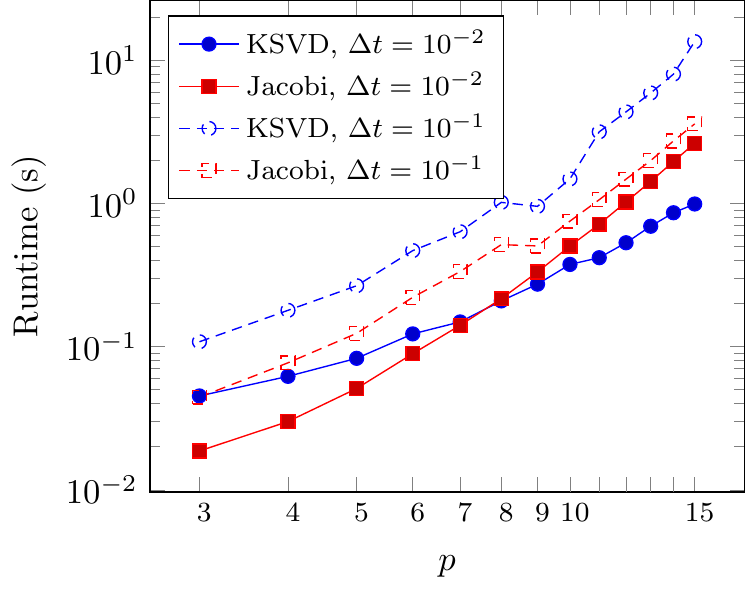}
        \caption{Total wall-clock time required per backward Euler step}
        \label{fig:ev-total-time}
    \end{subfigure}
    \caption{Runtime performance comparison of Kronecker-product preconditioner
             with exact block Jacobi preconditioner for 2D Euler equations.}
\end{figure}

In this section we compare the runtime performance of the Kronecker-product 
preconditioner with the exact block Jacobi preconditioner. Although we have 
observed that for large time steps $\Delta t$ or polynomial degrees $p$, the 
KSVD preconditioner requires more iterations to converge, it is also possible 
to compute and apply this preconditioner much more efficiently. Here, we compare
the wall-clock time required to compute and apply the preconditioner, according 
to Algorithms \ref{alg:2d-form} and \ref{alg:2d-apply}, respectively, for 
$p = 3, 4, \ldots, 15$. The block Jacobi preconditioner is computed by first 
assembling the diagonal block of the Jacobian matrix using the sum-factorized 
form of expression \eqref{eq:jac-entries}, and then computing its LU 
factorization. The wall-clock times for these operations are shown in Figure 
\ref{fig:ev-timings}. We remark that we observe the expected asymptotic 
computational complexities for each of these operations, where forming the 
Jacobi preconditioner requires $\mathcal{O}(p^6)$ operations, and applying 
the Jacobi preconditioner requires $\mathcal{O}(p^4)$ operations. Both forming 
and applying the approximate Kronecker-product preconditioner can be done in 
$\mathcal{O}(p^3)$ time. The total runtime observed per backward Euler step 
for $\Delta t = 0.1$ and $\Delta t = 0.01$ is shown in Figure 
\ref{fig:ev-total-time}. We see that for $p \geq 7$ and $\Delta t = 0.01$, the 
Kronecker-product preconditioner results in overall faster runtime, while for 
$\Delta t = 0.1$, because of the large number of iterations required per solve, 
the Jacobi preconditioner results in overall faster performance.

\subsection{2D Kelvin-Helmholtz instability}
For a more sophisticated test case, we consider a two-dimensional
Kelvin-Helmholtz instability. This important fluid instability occurs in shear
flows of fluids with different densities. The domain is taken to be the periodic
unit square $[0,1]^2$.
\begin{figure}[t!]
    \centering
    \begin{subfigure}[t]{0.5\textwidth}
        \centering
        \includegraphics[width=2.9in]{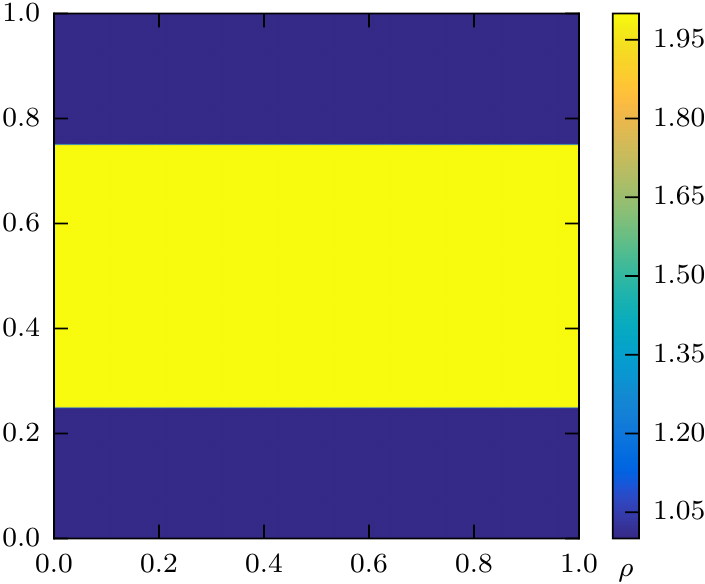}
        \caption{Density}
    \end{subfigure}%
    ~ 
    \begin{subfigure}[t]{0.5\textwidth}
        \centering
        \includegraphics[width=3.0in]{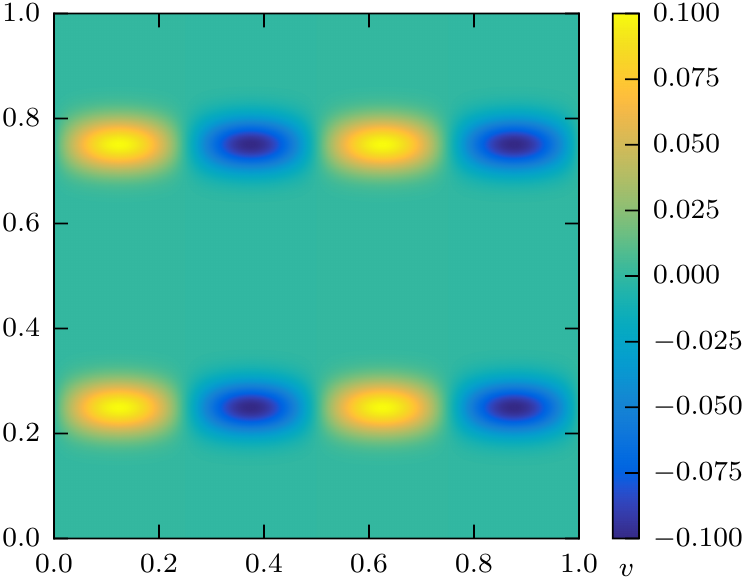}
        \caption{Vertical velocity}
        \label{fig:KH-IC-v}
    \end{subfigure}
    \caption{Initial conditions for the Kelvin-Helmholtz instability.}
    \label{fig:KH-IC}
\end{figure}
We define the function
\begin{equation}
f(x) = \frac{1}{4}(\erf(\alpha(x - 0.25)) + 1)(\erf(\alpha(0.75-x)) + 1),
\end{equation}
where $\alpha = 100$, as a smooth approximation to the discontinuous
characteristic function
\begin{equation}
  \chi(x) = \begin{cases}
  1, &\qquad 0.25 \leq x \leq 0.75 \\
  0, &\qquad \text{otherwise.}
  \end{cases}
\end{equation}
Following \cite{Springel:2010}, we define the initial conditions by
\begin{equation}
  \rho(x,y) = f(y) + 1, \qquad u(x,y) = f(y) - 1/2, 
  \qquad p(x,y) = 2.5,
\end{equation}
where the vertical velocity is given by
\begin{equation}
  v(x,y) = \frac{1}{10}\sin(4 \pi x)
           \left(\exp\left(-\frac{(y-0.25)^2}{2\sigma^2}\right)
               + \exp\left(-\frac{(y-0.75)^2}{2\sigma^2}\right)\right).
\end{equation}
Thus, the fluid density is equal to 2 inside the strip $y \in [0.25, 0.75]$, 
and 1 outside the strip. The fluid is moving to the right with horizontal 
velocity 0.5 inside the strip, and is moving to the left with equal speed 
outside of the strip. A small perturbation in the vertical velocity, localized 
around the discontinuity, determines the large-scale behavior of the 
instability. The initial conditions are shown in Figure \ref{fig:KH-IC}.

\begin{figure}[t!]
  \centering
    \begin{subfigure}[t]{0.5\textwidth}
        \centering
        \includegraphics[width=2.4in]{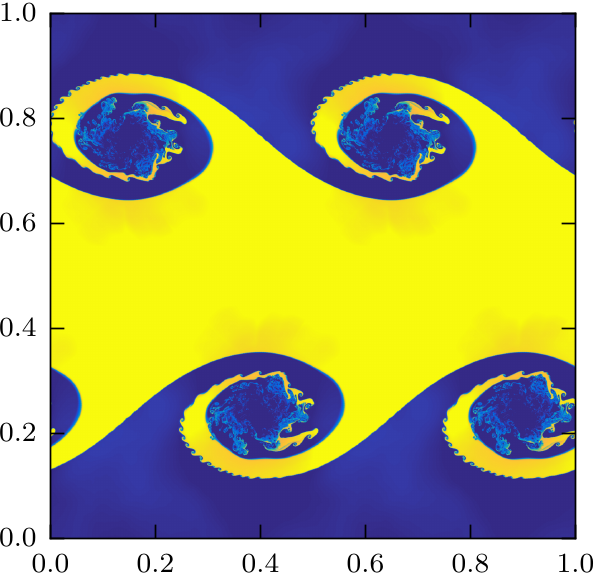}
    \end{subfigure}%
    ~ 
    \begin{subfigure}[t]{0.5\textwidth}
        \centering
        \includegraphics[width=2.5in]{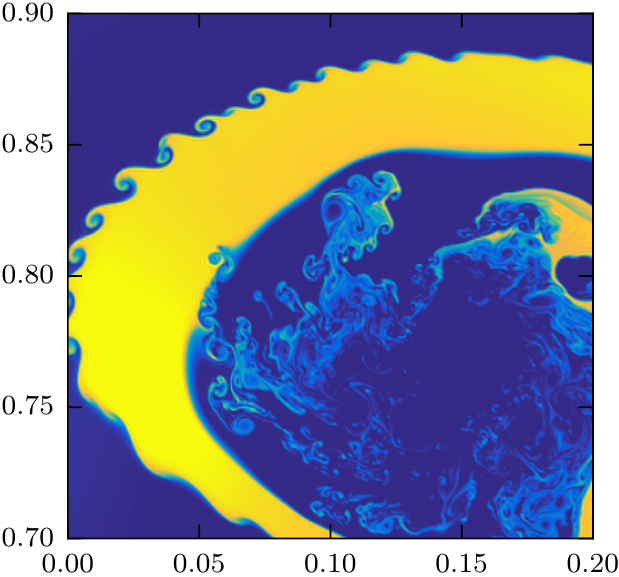}
    \end{subfigure}
  \caption{Solution (density) of Kelvin-Helmholtz instability at 
           $t = 1.5~\rm s$, with zoomed-in subregion to show small-scale 
           features.}
  \label{fig:KH-final}
\end{figure}

We use a $128 \times 128$ cartesian grid, with polynomial bases of degree 
3, 7, and 10. For 10th degree polynomials, the total number of degrees of 
freedom is 7,929,856. The Euler equations are integrated for $1.5~\rm s$ using a
fourth-order explicit method with a time step of $\Delta t = 2.5 \times 10^{-5}$
on the NERSC Edison supercomputer, using 480 cores. At this point, the solution
has developed sophisticated large- and small-scale features, as shown in Figure
\ref{fig:KH-final}.

We then linearize the Euler equations around this solution
in order to test the preconditioner performance. Because of the varied scale of
the features in this solution, we believe that the resulting linearization is a
representative of the DG systems we are interested in solving. Using this
solution, we then solve one backward Euler step using both the Jacobi and
approximate Kronecker-product preconditioners. For the implicit solve, we choose
a range of time steps, from the explicit step size of $\Delta t = 2.5 \times 
10^{-5}$, to a larger step of $\Delta t = 10^{-3}$. The number of 
GMRES iterations required per linear solve are show in Table 
\ref{tab:KH-results}. We observe that for the explicit-scale time step, the 
exact block Jacobi and approximation Kronecker-product preconditioner exhibit 
very similar performance for all choices of $p$. For the largest time step, 
$\Delta t = 10^{-3}$, the Kronecker-product preconditioner required about 
twice as many iterations for $p=3$, three times as many for $p=7$, and four 
times as many iterations for $p=10$.

\begin{table}[h!]
\caption{Number of GMRES iterations for Jacobi (\textit{J}) and KSVD 
         (\textit{K}) preconditioners, Euler equations for 2D Kelvin-Helmholtz
         instability.}
\label{tab:KH-results}
\centering
\begin{subtable}{0.3\linewidth}
\centering
\caption*{$p=3$}
\begin{tabular}{lll}
\toprule
$\Delta t$ & \textit{J} & \textit{K} \\
\midrule
$2.5 \times 10^{-5}$ & 4 & 4 \\
$5.0 \times 10^{-5}$ & 5 & 5 \\
$1.0 \times 10^{-4}$ & 6 & 7 \\
$2.5 \times 10^{-4}$ & 8 & 10 \\
$5.0 \times 10^{-4}$ & 10 & 14 \\
$1.0 \times 10^{-3}$ & 13 & 22 \\
\bottomrule
\end{tabular}
\end{subtable}
\begin{subtable}{0.3\linewidth}
\centering
\caption*{$p=7$}
\begin{tabular}{lll}
\toprule
$\Delta t$ & \textit{J} & \textit{K} \\
\midrule
$2.5 \times 10^{-5}$ & 5 & 6 \\
$5.0 \times 10^{-5}$ & 6 & 8 \\
$1.0 \times 10^{-4}$ & 8 & 12 \\
$2.5 \times 10^{-4}$ & 12 & 20 \\
$5.0 \times 10^{-4}$ & 15 & 35 \\
$1.0 \times 10^{-3}$ & 21 & 62 \\
\bottomrule
\end{tabular}
\end{subtable}
\begin{subtable}{0.3\linewidth}
\centering
\caption*{$p=10$}
\begin{tabular}{lll}
\toprule
$\Delta t$ & \textit{J} & \textit{K} \\
\midrule
$2.5 \times 10^{-5}$ & 6 & 8 \\
$5.0 \times 10^{-5}$ & 8 & 11 \\
$1.0 \times 10^{-4}$ & 10 & 16 \\
$2.5 \times 10^{-4}$ & 14 & 31 \\
$5.0 \times 10^{-4}$ & 19 & 55 \\
$1.0 \times 10^{-3}$ & 24 & 106 \\
\bottomrule
\end{tabular}
\end{subtable}
\end{table}

{
\subsection{2D NACA Airfoil}

In this test case, we consider the viscous flow over a NACA 0012 airfoil with
angle of attack $30^\circ$ at Reynolds number 16000, with Mach number $M_0 =
0.2$. We take the domain to be a disk of radius 10, centered at (0,0). The
leading edge of the airfoil is placed at the origin. A no-slip wall condition is
enforced at the surface of the airfoil, and far-field conditions are enforced at
all other domain boundaries. The far-field velocity is set to unity in the
freestream direction. The domain is discretized using an unstructured
quadrilateral mesh, refined in the vicinity of the wing and in its wake.
Isoparametric mappings are used to curve the elements on the airfoil surface.
This flow is characterized by the thin boundary layer that develops on the
airfoil. In order to resolve this boundary layer, we introduce stretched,
anisotropic boundary-layer elements at the surface of the airfoil. These small
elements result in a CFL condition that requires the use of very
small time steps when using an explicit time integration method. The mesh and
density contours are shown in Figure \ref{fig:naca-mesh}.

\begin{figure}[t]
\centering
\begin{subfigure}[t]{\textwidth}
\centering
\includegraphics[width=5in]{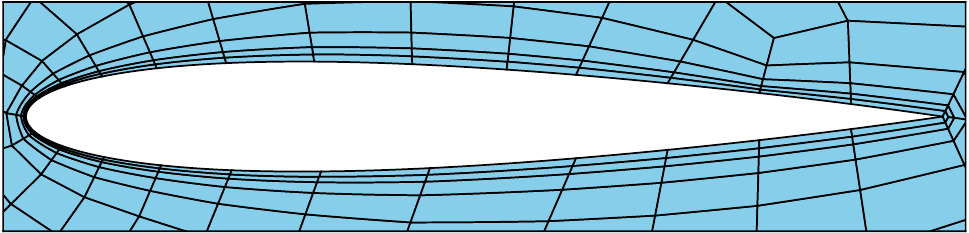}
\end{subfigure}
\\[0.1in]
\begin{subfigure}[t]{0.49\textwidth} \centering
\includegraphics[width=2.8in]{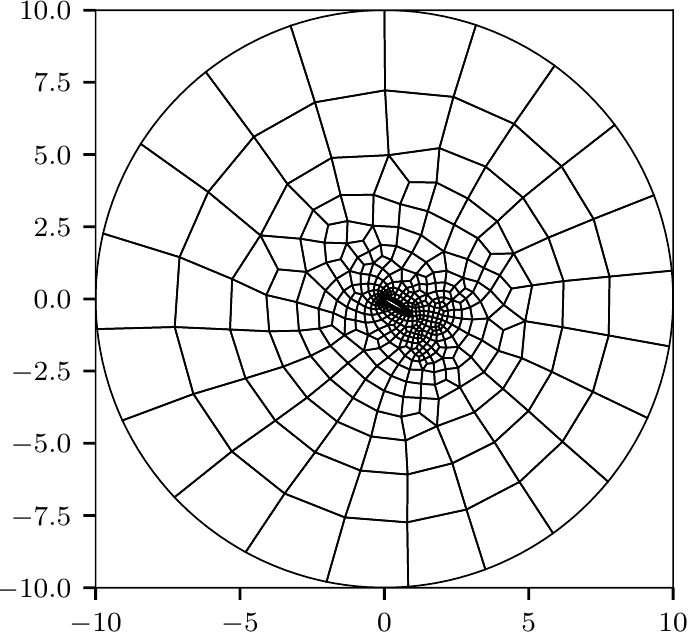}
\end{subfigure}
\begin{subfigure}[t]{0.49\textwidth} \centering
\includegraphics[width=3.1in]{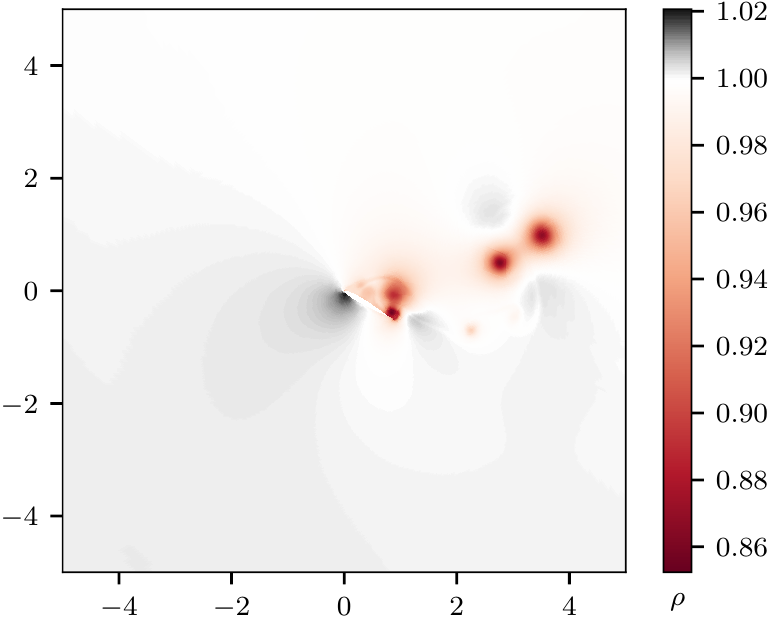}
\end{subfigure}
\caption{NACA 0012 mesh, with zoom-in around the surface of the airfoil showing
         anisotropic boundary-layer elements (top), and solution (density).}
\label{fig:naca-mesh}
\end{figure}

This test case differs from the preceding two test cases because instead of the 
Euler equations we solve the compressible Navier-Stokes equations,
\begin{gather}
    \label{eq:ns-1}
    \frac{\partial\rho}{\partial t}+\frac{\partial}{\partial x_j}(\rho u_j)=0\\
    \label{eq:ns-2}
    \frac{\partial}{\partial t}(\rho u_i) + \frac{\partial}{\partial x_j}
        (\rho u_i u_j) + \frac{\partial p}{\partial x_i1}
        = \frac{\partial \tau_{ij}}{\partial x_j}
        \qquad \text{for $i=1,2,3,$}\\
    \label{eq:ns-3}
    \frac{\partial}{\partial t}(\rho E) + \frac{\partial}{\partial x_j}
        \left(u_j(\rho E + p) \right) = -\frac{\partial q_j}{\partial x_j}
        + \frac{\partial}{\partial x_j} (u_j \tau_{ij}),
\end{gather}
The viscous stress tensor and heat flux are given by
\begin{equation}
\tau_{ij} = \mu\left( \frac{\partial u_i}{\partial x_j} +
                      \frac{\partial u_j}{\partial x_i} -
                      \frac{2}{3} \frac{\partial u_k}{\partial x_k} \delta_{ij}
                      \right)
\qquad\text{and}\qquad
q_j = - \frac{\mu}{\rm Pr} \frac{\partial}{\partial x_j} \left(
E + \frac{p}{\rho} - \frac{1}{2} u_k u_k \right).
\end{equation}
Here $\mu$ is the coefficient of viscosity, and ${\rm Pr}$ is the Prandtl 
number, which we assume to be constant. We discretize the second-order terms 
using the local discontinuous Galerkin method \cite{Cockburn:1998:LDG}, which 
introduces certain lifting operators into the primal form of the discretization
\cite{Arnold:2002}. These lifting operators do not readily fit into the 
tensor-produce framework described above, and thus we apply the approximate 
Kronecker-product preconditioner to only the inviscid component of the flux 
function. Since the flow is convection-dominated away from the airfoil, we 
believe that this provides an acceptable approximation, although properly 
incorporating the viscous terms into the preconditioner is an area of ongoing 
research.

We integrate the equations until $t=2.5$ in order to obtain a
representative initial condition about which to linearize the equations. We
consider polynomial degrees $p = 1, 3, 7, 10, 15$, and compare the efficiency of
explicit and implicit time integration methods, using both the Kronecker-product
preconditioner and exact block Jacobi. In order to make this comparison, we
determine experimentally the largest explicit timestep for which the system is
stable. Then, we measure the wall-clock time required to integrate the system
from $t=2.5$ until $t=3.5$. Similarly, for the implicit methods,
we experimentally choose an appropriate $\Delta t$ measure the wall-clock time
required to advance the simulation until $t=3.5$ using a
three-stage, third-order, $L$-stable DIRK scheme. We present these results in
Table \ref{tab:naca-times}. We note that for $p > 3$, the Kronecker-product
preconditioner results in the shortest runtimes. For large $p$, the exact block
Jacobi preconditioner becomes impractical due to the large $p$-dependence of the
computational complexity. For $p=15$, the high degree polynomials considered for
this test case, the Kronecker preconditioner resulted in runtimes that were
about a factor of two faster than explicit, and a factor of ten faster than
exact block Jacobi.

\begin{table}[H] {
\setlength{\tabcolsep}{12pt}
\begin{tabular}{llllll}
\toprule
&&&\multicolumn{3}{c}{Runtime (s)} \\
\cmidrule{4-6}
$p$ & Explicit $\Delta t$ & Implicit $\Delta t$ 
    & \multicolumn{1}{c}{RK4} & \multicolumn{1}{c}{K} & \multicolumn{1}{c}{J} \\
\midrule
1 & $1.0 \times 10^{-4}$ & $2.0 \times 10^{-1}$ & $8.864 \times 10^1$ 
                                  & $5.225 \times 10^1$ & $5.854 \times 10^1$ \\
3 & $2.5 \times 10^{-5}$ & $5.0 \times 10^{-2}$ & $6.381 \times 10^2$ 
                                  & $8.406 \times 10^2$ & $4.626 \times 10^2$ \\
7 & $1.0 \times 10^{-6}$ & $2.0 \times 10^{-3}$ & $5.441 \times 10^4$ 
                                  & $3.057 \times 10^4$ & $6.765 \times 10^4$ \\
10 & $5.0 \times 10^{-7}$ & $2.5 \times 10^{-4}$ & $2.218 \times 10^5$ 
                                  & $1.884 \times 10^5$ & $1.949 \times 10^6$ \\
15 & $1.0 \times 10^{-7}$ & $2.0 \times 10^{-4}$ & $2.372 \times 10^6$ 
                                  & $1.176 \times 10^6$ & $1.223 \times 10^7$ \\
\bottomrule
\end{tabular} }
\caption{Runtime results for NACA test-case, comparing explicit (fourth-order 
         Runge-Kutta) with implicit (three-stage DIRK), using Kronecker and 
         exact block Jacobi preconditioners. Runtime in seconds per simulation 
         second is presented.}
\label{tab:naca-times}
\end{table}

}

\subsection{3D periodic Euler}\label{sec:periodic-euler}
\begin{figure}[t!]
 \centering
 \includegraphics[width=3.0in]{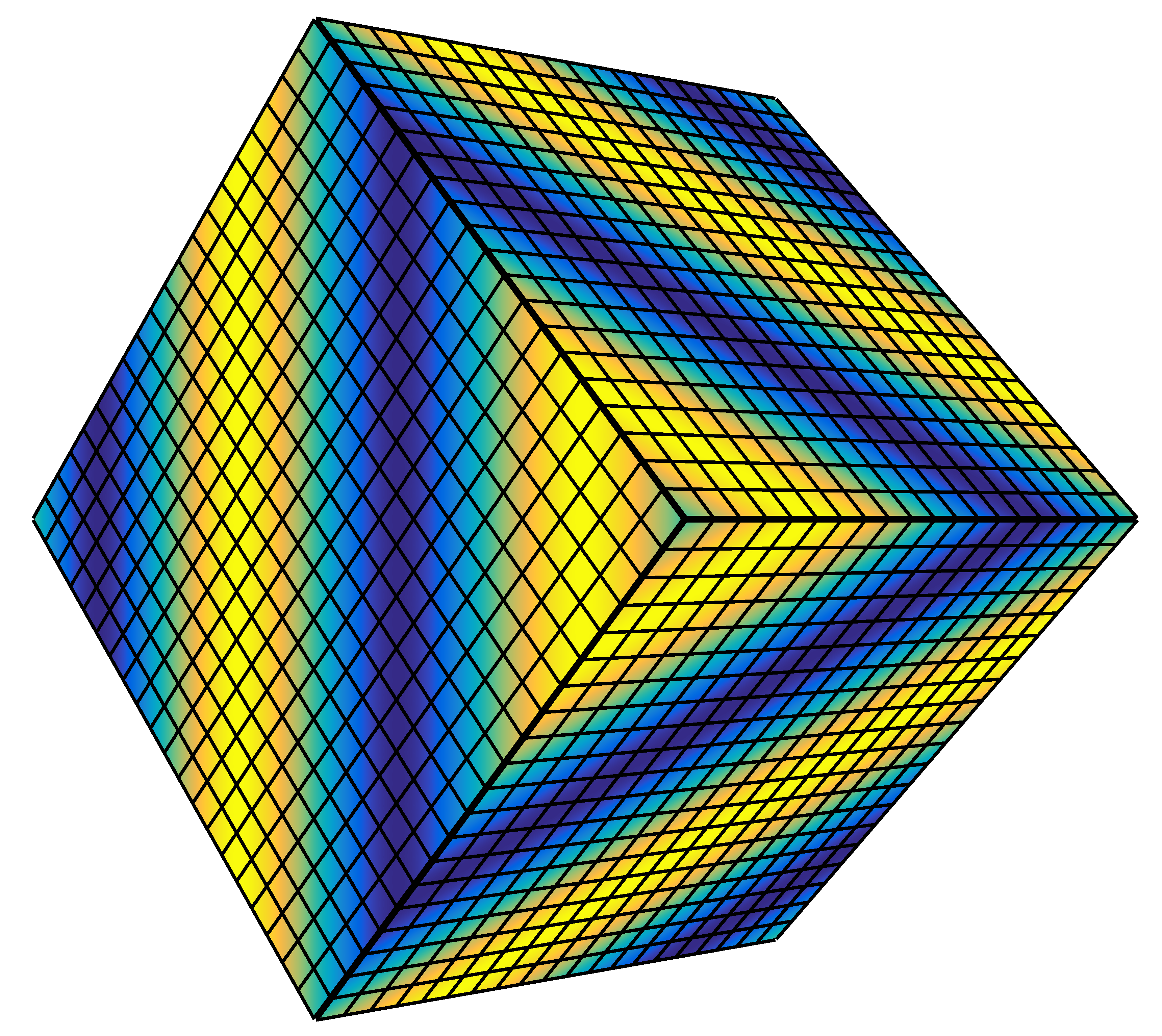}
 \caption{Initial conditions (density) for smooth 3D Euler test case.}
 \label{fig:3d-periodic-ic}
\end{figure}

In this example, we provide a test case for the three-dimensional 
preconditioner. We consider the cube $[0,2]^3$ with periodic boundary 
conditions, and use the initial conditions from \cite{Jiang:1996}, given by
\begin{align}
 &\rho = 1 + 0.2 \sin(\pi(x + y + z)), \\
 &u = 1, \quad v = -1/2, \quad w = 1, \\
 &p = 1.
\end{align}
The exact solution to the Euler equations is known analytically in this case.
Velocity and pressure remain constant in time, and the density at time $t$ is 
given by
\begin{equation}
 \rho = 1 + 0.2\sin(\pi(x + y + z - t(u + v + w))).
\end{equation}
The initial conditions are shown in Figure \ref{fig:3d-periodic-ic}. The mesh is
taken to be a regular $6 \times 6 \times 6$ hexahedral grid, and we consider
polynomial degrees of $p = 1, 2, \ldots, 12$. We choose a representative time
step of $\Delta t = 2.5 \times 10^{-3}$. This time step can be used for explicit
methods with low-degree polynomials, but because of the $p$-dependency of the
CFL condition, we observe that for $p \geq 9$, explicit methods become unstable,
motivating the use of implicit methods. In order to compare the efficiency of
the approximate Kronecker-product preconditioner, we compute one backward Euler
step and compare the number of GMRES iterations required per linear solve using
the exact block Jacobi preconditioner and the KSVD preconditioner.

\begin{figure}[h!]
\begin{floatrow}
\capbtabbox{%
 \centering
 \begin{tabular}{l|llll}
 \toprule
 $p$ & $J_{\rm full}$ & $K_{\rm full}$ & $J_{\rm small}$ & $K_{\rm small}$ \\
 \midrule
 1 & 4 & 4 & 5 & 5 \\
 2 & 4 & 5 & 6 & 6 \\
 3 & 5 & 5 & 7 & 7 \\
 4 & 5 & 5 & 8 & 8 \\
 5 & 5 & 6 & 9 & 10 \\
 6 & 5 & 6 & 11 & 12 \\
 7 & 5 & 6 & 12 & 13 \\
 8 & 5 & 7 & 15 & 16 \\
 9 & 5 & 7 & 16 & 17 \\
 10 & 5 & 8 & 19 & 21 \\
 11 & 5 & 8 & 20 & 22 \\
 12 & - & 9 & 23 & 26 \\
 \bottomrule
 \end{tabular}
 \vspace{0.15in}
}{%
 \caption{GMRES iterations required per linear solve for three-dimensional 
          periodic Euler test case.}%
 \label{tab:periodic-results}%
} \hspace{0.4in}
\ffigbox{%
 \centering
 \includegraphics[width=3in]{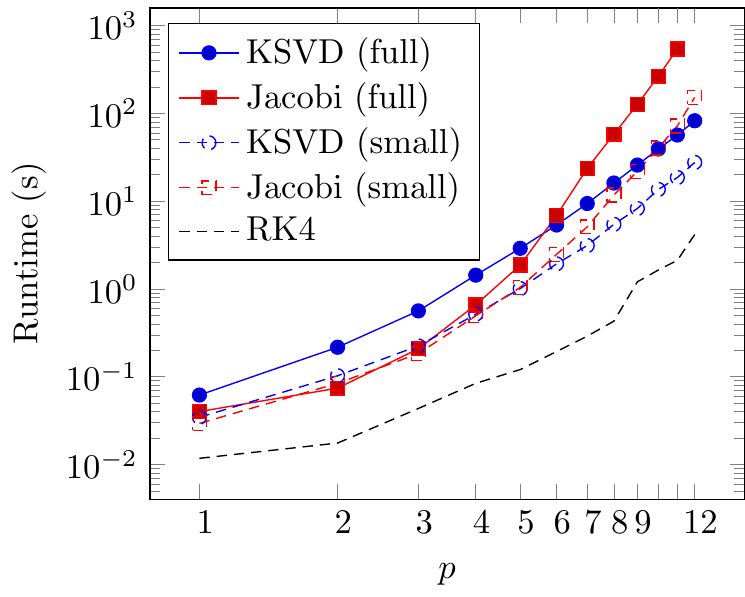}%
}{%
 \caption{Wall-clock time required per backward Euler solve for 
          three-dimensional periodic Euler test case.
          For reference, wall-clock
          time required for explicit RK4 is shown.}%
 \label{fig:periodic-timings}%
}
\end{floatrow}
\end{figure}

In this test case, we also consider two variations each of these
preconditioners. Since the solution to the Euler equations consists of five
components, we can consider the diagonal blocks of the Jacobian to either be
large $5 (p+1)^3 \times 5 (p+1)^3$ blocks coupling all solution components, or
as smaller $(p+1)^3 \times (p+1)^3$ blocks, which do not couple the solution
components. We expect that using the smaller blocks will require more GMRES
iterations per linear solve, because each block captures less information. The
larger blocks, on the other hand, are much more computationally expensive to
compute. We present the number of iterations required to converge for each 
of the preconditioners in Table \ref{tab:periodic-results}, where $J$ and $K$ 
stand for the block Jacobi and Kronecker-product preconditioners, respectively,
and the subscripts ``full'' and ``small'' refer to the block size used. We 
observe that for small polynomial degree $p$, the KSVD preconditioner results 
in close-to-identical number of iterations, when compared with exact block 
Jacobi. For this test case, for $p$ closer to 12, we observe that the number of 
iterations grows faster for the KSVD preconditioner than for the block Jacobi 
preconditioner, but the difference remains relatively small. Additionally, as 
expected, the small-block preconditioner requires a greater number of iterations
to converge when compared with the full-block preconditioner. We note that due 
to the very large memory requirements, the full Jacobi preconditioner did not 
complete for $p=12$.

We additionally measure the wall-clock time required per backward Euler step for
each of the preconditioners, and present the results in Figure
\ref{fig:periodic-timings}. For reference, we also include the 
wall-clock time required to integrate the system of equations for an equivalent
time using the explicit RK4 method with a stable time step. This problem is 
well-suited for explicit solvers, and thus the RK4 method is more efficient than
the implicit methods.
From these measurements, it is possible to draw
several conclusions. Firstly, as $p$ grows, it becomes possible to
observe the $\mathcal{O}(p^9)$ complexity of the exact block Jacobi
preconditioner, which becomes prohibitively expensive. This is in contrast to
the Kronecker-product preconditioner, whose $\mathcal{O}(p^5)$ complexity
results in reasonable runtimes for all $p$ considered. The full-block KSVD
preconditioner results in faster runtimes starting at about $p=5$, and the
small-block preconditioner at about $p=4$. We also see that, despite the larger
number of iterations required, the small-block preconditioner results in faster
overall runtime.

{
\subsection{Compressible Taylor-Green Vortex at ${\rm Re} = 1600$}
The direct numerical simulation of the Taylor-Green vortex at ${\rm Re} = 1600$
is a benchmark problem from the first International Workshop on High-Order CFD
Methods \cite{Wang2013}. This three-dimensional problem has been often used to 
study the performance of high-order methods, and DG methods in particular, on 
transitional flows. This problem provides a useful test case because of the 
availability of fully-resolved reference data
\cite{Shu2005,van2011comparison,Carton2014assessment,DeBonis2013,Chapelier2012}.
The domain is taken to be the periodic cube $[-\pi,\pi]^3$. The initial 
conditions are given by
\begin{align}
  u(x,y,z) &= u_0 \sin(x) \cos(y) \cos(z) \\
  v(x,y,z) &= -u_0 \cos(x) \sin(y) \cos(z) \\
  w(x,y,z) &= 0 \\
  p(x,y,z) &= p_0 + \frac{\rho_0 u_0^2}{16}\left(
    \cos(2x) + \cos(2y) \right)\left(\cos(2z) + 2\right),
\end{align}
where the parameters are given by $\gamma = 1.4$, ${\rm Pr} = 0.71$, $u_0 = 1$,
$\rho_0 = 1$, with Mach number $M_0 = u_0/c_0 = 0.10$, where $c_0$ is the speed
of sound computed in accordance with the pressure $p_0$. The initial density
distribution is then given by $\rho = p \rho_0 / p_0$. The characteristic
convective time is given by $t_{\rm c}=1$, and the final time is $t_{\rm f} = 20
t_{\rm c}$. The geometry is discretized using regular hexahedral grids of size 
$8^3, 16^3, 32^3,$ and $42^3$, with polynomial degrees $p = 3, 4, 7,$ and $15$.

\begin{figure}[t]
\centering
\includegraphics[width=5.5in]{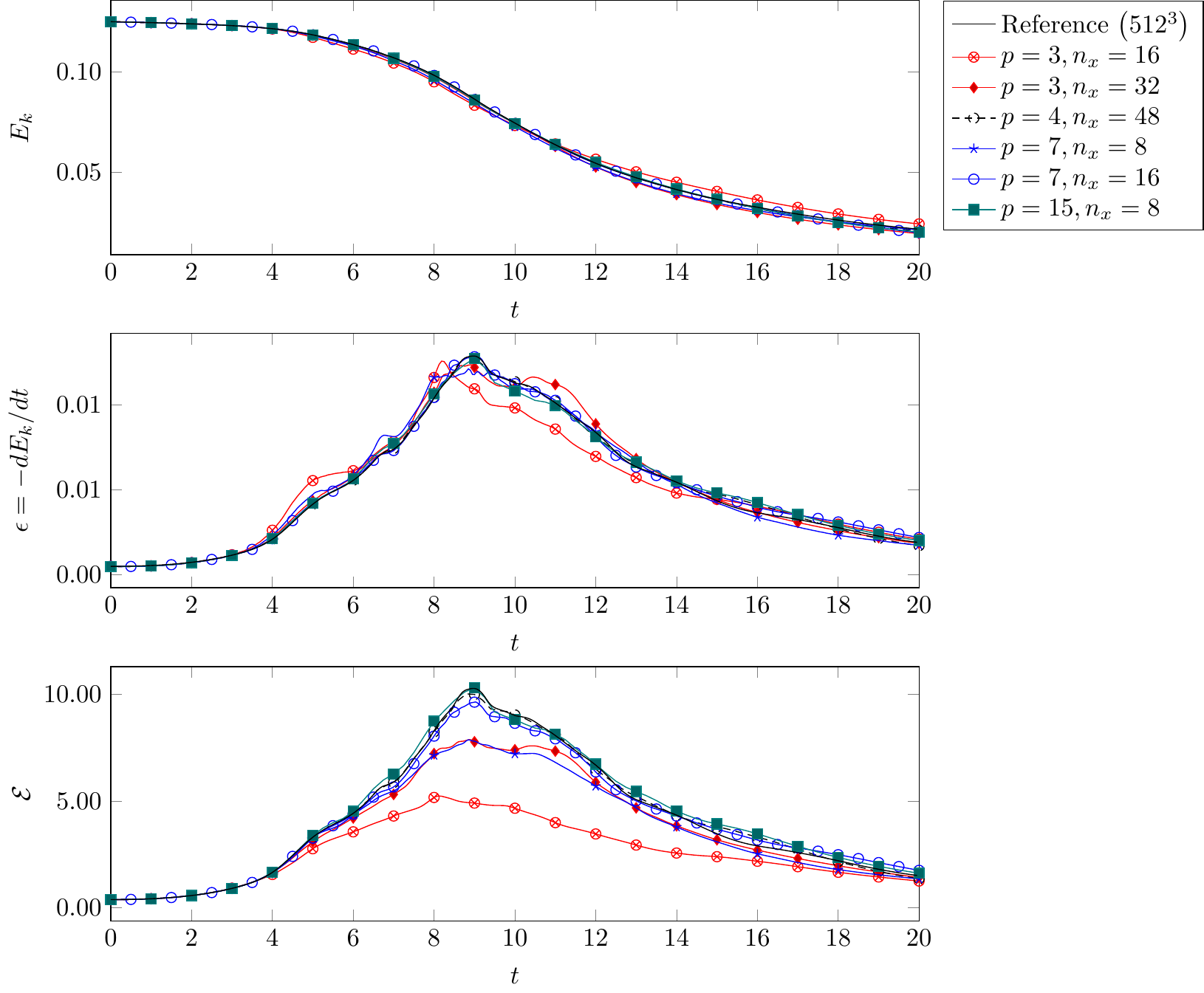}
\caption{Time evolution of kinetic energy $E_k$, kinetic energy dissipation rate
         (KEDR) $\epsilon$, and enstrophy $\mathcal{E}$, for Taylor-Green test 
         case. Comparison of various DG configurations with reference 
         pseudo-spectral solution.}
\label{fig:tg-diagnostics}
\end{figure}

The Taylor-Green vortex provides a strong motivation for the use of very high 
polynomial degrees. In Figure \ref{fig:tg-diagnostics}, we show the 
time-evolution of the diagnostic quantities of mean energy, kinetic energy 
dissipation rate, and enstrophy,
\begin{align}
  E_k(t) &= \frac{1}{\rho_0|\Omega|}
    \int_\Omega \rho \frac{\bm{u}\cdot \bm{u}}{2}~dx, \\
  \epsilon(t) &= -\frac{dE_k}{dt}(t), \\
  \mathcal{E}(t) &= \frac{1}{\rho_0|\Omega|}
    \int_\Omega \rho \frac{\omega \cdot \omega}{2}~dx.
\end{align}
For each grid configuration, we compare the results with a fully-resolved 
pseudo-spectral reference solution. We notice that the relatively low-order 
solutions with $p=3$ severely underpredict the peak enstrophy. However, with 
equal numbers of numerical degrees of freedom, the higher-order $p=7$ and $p=15$
solutions much more closely match the reference data. For example, the $n_x = 8,
p = 15$ discretization results in much better agreement than the $n_x=32, p=3$ 
case, despite an equal number of degrees of freedom, motivating the use of 
very high polynomial degrees.

We now examine the efficiency of the approximate tensor-product preconditioner
compared with exact block Jacobi for each of the grid configurations shown in
Figure \ref{fig:tg-diagnostics}. For each configuration, we choose a range of
timesteps, ranging from $\Delta t = 10^{-4}$ to $\Delta t = 1.6 \times 10^{-3}$
by factors of two. We measure the average number of GMRES iterations per linear
solve, and list the results in Table \ref{tab:tg-iters}. In the case of the
$p=15, n_x=8$, the exact block Jacobi preconditioner did not complete because of
excessive runtime and memory requirements. These iteration counts indicate that
the number of GMRES iterations per linear solve increases as the timestep
increases, and that this dependence is sublinear. For almost all cases, the
Kronecker-product preconditioner resulted in very similar to iteration counts
when compared with exact block Jacobi. However, due to the decreased
computational complexity and decreased memory requirements, the KSVD
preconditioner can be used to obtain similar GMRES convergence at a highly
decreased computational cost.

\begin{figure}[b!]
    \centering
    \includegraphics[width=3in]{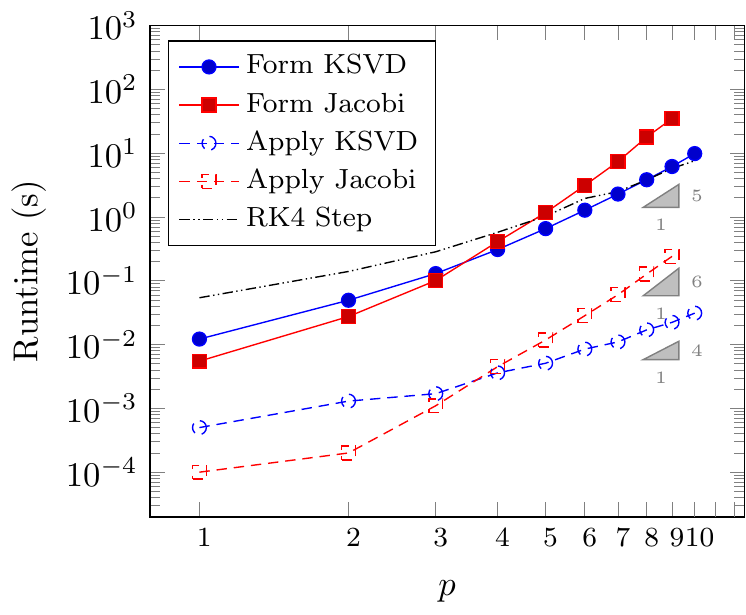}
    \caption{
             Wall-clock time required to form (solid lines) and 
             apply (dashed lines) the exact block Jacobi and approximate 
             Kronecker-product preconditioners. For reference, also shown is 
             time per RK4 step.}
    \label{fig:tg-form-apply}
\end{figure}

In Figure \ref{fig:tg-form-apply}, we show the wall-clock times required to form
and apply the approximate Kronecker and exact block Jacobi preconditioners with
$n_x = 8$ and $p = 1, 2, \ldots, 10.$ Due to excessive memory requirements, the
exact Jacobi preconditioner was not computed for $p=10$. The improved
computational complexities for the KSVD preconditioner for both operations are
apparent. Forming the Kronecker-product preconditioner requires
$\mathcal{O}(p^5)$ operations and applying the preconditioner requires
$\mathcal{O}(p^4)$ operations, as opposed to $\mathcal{O}(p^9)$ and
$\mathcal{O}(p^6)$, respectively, for the block Jacobi preconditioner. For
comparison, the wall-clock time required to perform one explicit RK4 step is
also shown. By taking advantage of the tensor-product structure, as described in
Section \ref{sec:explicit}, the computational complexity of performing an
explicit step scales as $\mathcal{O}(p^4)$. However, the choice of stable time
step is severely restricted as $p$ grows, requiring a large number of steps to
be taken. This problem is structurally quite similar to that of Section
\ref{sec:periodic-euler}, since the geometry for both problems is a regular
cartesian grid, and the equations only differ in the presence of viscosity.
Therefore, the performance characteristics for the solvers are quite similar to
those shown in Figure \ref{fig:periodic-timings}. For example, we observe that
for $p=7, n_x = 8$, the explicit time integration results in an overall runtime
that is about one-fifth of the runtime for implicit time integration with the
KSVD preconditioner, and one-seventh of the runtime for implicit time
integration with the block Jacobi preconditioner.

\begin{table}[h!]
\centering
\caption{Average number of GMRES iterations per linear solve for Taylor-Green
         test case.}
\label{tab:tg-iters}
\begin{subtable}{0.3\linewidth}
\centering
\caption{$p=3, n_x = 16$}
\begin{tabular}{lll}
\toprule
$\Delta t$ & Jacobi & KSVD \\
\midrule
$1\times 10^{-4}$   & 4  & 4\\
$2\times 10^{-4}$   & 5  & 7\\
$4\times 10^{-4}$   & 6  & 7\\
$8\times 10^{-4}$   & 8  & 9\\
$1.6\times 10^{-3}$ & 11 & 12\\
\bottomrule
\end{tabular}
\end{subtable} \hspace{0.5cm}
\begin{subtable}{0.3\linewidth}
\centering
\caption{$p=3, n_x = 32$}
\begin{tabular}{lll}
\toprule
$\Delta t$ & Jacobi & KSVD \\
\midrule
$1\times 10^{-4}$   & 5  & 5\\
$2\times 10^{-4}$   & 6  & 6\\
$4\times 10^{-4}$   & 7  & 8\\
$8\times 10^{-4}$   & 10 & 12\\
$1.6\times 10^{-3}$ & 16 & 18\\
\bottomrule
\end{tabular}
\end{subtable} \hspace{0.5cm}
\begin{subtable}{0.3\linewidth}
\centering
\caption{$p=4, n_x = 48$}
\begin{tabular}{lll}
\toprule
$\Delta t$ & Jacobi & KSVD \\
\midrule
$1\times 10^{-4}$   & 7  & 7 \\
$2\times 10^{-4}$   & 9  & 10\\
$4\times 10^{-4}$   & 12 & 16\\
$8\times 10^{-4}$   & 19 & 27\\
$1.6\times 10^{-3}$ & 30 & 40\\
\bottomrule
\end{tabular}
\end{subtable}

\vspace{0.25in}

\begin{subtable}{0.3\linewidth}
\centering
\caption{$p=7, n_x = 8$}
\begin{tabular}{lll}
\toprule
$\Delta t$ & Jacobi & KSVD \\
\midrule
$1\times 10^{-4}$   & 5  & 6\\
$2\times 10^{-4}$   & 6  & 6\\
$4\times 10^{-4}$   & 8  & 8\\
$8\times 10^{-4}$   & 10 & 11\\
$1.6\times 10^{-3}$ & 15 & 16\\
\bottomrule
\end{tabular}
\end{subtable} \hspace{0.5cm}
\begin{subtable}{0.3\linewidth}
\centering
\caption{$p=7, n_x = 16$}
\begin{tabular}{lll}
\toprule
$\Delta t$ & Jacobi & KSVD \\
\midrule
$1\times 10^{-4}$   & 6  & 6\\
$2\times 10^{-4}$   & 7  & 9\\
$4\times 10^{-4}$   & 10 & 13\\
$8\times 10^{-4}$   & 15 & 17\\
$1.6\times 10^{-3}$ & 25 & 29\\
\bottomrule
\end{tabular}
\end{subtable} \hspace{0.5cm}
\begin{subtable}{0.3\linewidth}
\centering
\caption{$p=15, n_x = 8$}
\begin{tabular}{lll}
\toprule
$\Delta t$ & Jacobi & KSVD \\
\midrule
$1\times 10^{-4}$   & -- & 8\\
$2\times 10^{-4}$   & -- & 11\\
$4\times 10^{-4}$   & -- & 15\\
$8\times 10^{-4}$   & -- & 26\\
$1.6\times 10^{-3}$ & -- & 88\\
\bottomrule
\end{tabular}
\end{subtable}
\end{table}
}

\section{Conclusion and future work}
In this work we have developed new approximate tensor-product based
preconditioners for very high-order discontinuous Galerkin methods. These
preconditioners are computed using an algebraic singular value-based algorithm,
and compare favorably with the traditional block Jacobi preconditioner. The
computational complexity is reduced from $\mathcal{O}(p^6)$ to
$\mathcal{O}(p^3)$ in two spatial dimensions and $\mathcal{O}(p^9)$ to
$\mathcal{O}(p^5)$ in three spatial dimensions. Numerical results in two and
three dimensions for the advection and Euler equations, using polynomial degrees
up to $p=30$, confirm the expected computational complexities, and demonstrate
significant reductions in runtimes for certain test problems.

Future work for further improving the performance of this preconditioner include
systematic treatment of viscous fluxes and second-order terms, and fast 
inversion of sums of more than two Kronecker products allowing for 
treatment of off-diagonal blocks in the context of an ILU-based preconditioner. 
Also of interest is the investigation of the performance of the preconditioner 
when used as a smoother in $p$-multigrid solvers.

\section{Acknowledgments}
This research used resources of the National Energy Research Scientific
Computing Center, a DOE Office of Science User Facility supported by the
Office of Science of the U.S. Department of Energy under Contract
No.\ DE-AC02-05CH11231, and was supported by the AFOSR Computational
Mathematics program under grant number FA9550-15-1-0010. The first author was
supported by the Department of Defense through the National Defense Science \&
Engineering Graduate Fellowship Program and by the Natural Sciences and
Engineering Research Council of Canada.

\bibliographystyle{plain}
\bibliography{KSVD}

\begin{thebibliography}{10}

\bibitem{Alexander:1977dk}
Roger Alexander.
\newblock {Diagonally implicit Runge-Kutta methods for stiff O.D.E.'s}.
\newblock {\em SIAM Journal on Numerical Analysis}, 14(6):1006--1021, 1977.

\bibitem{Arnold:2002}
Douglas~N. Arnold, Franco Brezzi, Bernardo Cockburn, and L.~Donatella Marini.
\newblock Unified analysis of discontinuous {G}alerkin methods for elliptic
  problems.
\newblock {\em SIAM Journal on Numerical Analysis}, 39(5):1749--1779, 2002.

\bibitem{baggag1999parallel}
Abdalkader Baggag, Harold Atkins, and David Keyes.
\newblock Parallel implementation of the discontinuous galerkin method.
\newblock 1999.

\bibitem{Birken:2013}
Philipp Birken, Gregor Gassner, Mark Haas, and Claus-Dieter Munz.
\newblock Preconditioning for modal discontinuous {G}alerkin methods for
  unsteady {3D} {N}avier-{S}tokes equations.
\newblock {\em Journal of Computational Physics}, 240:20--35, may 2013.

\bibitem{Carton2014assessment}
C.~Carton~de Wiart, K.~Hillewaert, M.~Duponcheel, and G.~Winckelmans.
\newblock Assessment of a discontinuous {G}alerkin method for the simulation of
  vortical flows at high {R}eynolds number.
\newblock {\em International Journal for Numerical Methods in Fluids},
  74(7):469--493, 2014.

\bibitem{Chapelier2012}
Jean-Baptiste Chapelier, Marta De La~Llave Plata, and Florent Renac.
\newblock Inviscid and viscous simulations of the {T}aylor-{G}reen vortex flow
  using a modal discontinuous {G}alerkin approach.
\newblock In {\em 42nd {AIAA} Fluid Dynamics Conference and Exhibit}. American
  Institute of Aeronautics and Astronautics, June 2012.

\bibitem{Cockburn:1998:LDG}
Bernardo Cockburn and Chi-Wang Shu.
\newblock The local discontinuous {G}alerkin method for time-dependent
  convection-diffusion systems.
\newblock {\em SIAM Journal on Numerical Analysis}, 35(6):2440--2463, 1998.

\bibitem{Cockburn:1998:RKDG}
Bernardo Cockburn and Chi-Wang Shu.
\newblock The {R}unge-{K}utta discontinuous {G}alerkin method for conservation
  laws {V}: multidimensional systems.
\newblock {\em Journal of Computational Physics}, 141(2):199--224, 1998.

\bibitem{Crivellini2011}
A.~Crivellini and F.~Bassi.
\newblock An implicit matrix-free discontinuous {G}alerkin solver for viscous
  and turbulent aerodynamic simulations.
\newblock {\em Computers {\&} Fluids}, 50(1):81--93, nov 2011.

\bibitem{DeBonis2013}
James DeBonis.
\newblock Solutions of the {T}aylor-{G}reen vortex problem using
  high-resolution explicit finite difference methods.
\newblock In {\em 51st {AIAA} Aerospace Sciences Meeting including the New
  Horizons Forum and Aerospace Exposition}. American Institute of Aeronautics
  and Astronautics, January 2013.

\bibitem{Diosady:2011domain}
Laslo~T. Diosady.
\newblock {\em Domain Decomposition Preconditioners for Higher-Order
  Discontinuous Galerkin Discretizations}.
\newblock PhD thesis, Massachusetts Institute of Technology, 2011.

\bibitem{Diosady:2017}
Laslo~T. Diosady and Scott~M. Murman.
\newblock Tensor-product preconditioners for higher-order space–time
  discontinuous {G}alerkin methods.
\newblock {\em Journal of Computational Physics}, 330:296 -- 318, 2017.

\bibitem{EscobarVargas}
J.~A. Escobar-Vargas, P.~J. Diamessis, and C.~F.~Van Loan.
\newblock The numerical solution of the pressure {P}oisson equation for the
  incompressible {N}avier-{S}tokes equations using a quadrilateral spectral
  multidomain penalty method.
\newblock Preprint available at
  https://www.cs.cornell.edu/cv/ResearchPDF/Poisson.pdf, 2011.

\bibitem{Golub:1981}
Gene~H. Golub, Franklin~T. Luk, and Michael~L. Overton.
\newblock A block {L}anczos method for computing the singular values and
  corresponding singular vectors of a matrix.
\newblock {\em ACM Transactions on Mathematical Software (TOMS)},
  7(2):149--169, 1981.

\bibitem{Gopalakrishnan:2003}
J.~Gopalakrishnan and G.~Kanschat.
\newblock A multilevel discontinuous galerkin method.
\newblock {\em Numerische Mathematik}, 95(3):527--550, 2003.

\bibitem{gottlieb1991cfl}
David Gottlieb and Eitan Tadmor.
\newblock The {CFL} condition for spectral approximations to hyperbolic
  initial-boundary value problems.
\newblock {\em Mathematics of Computation}, 56(194):565--588, 1991.

\bibitem{Jiang:1996}
Guang-Shan Jiang and Chi-Wang Shu.
\newblock Efficient implementation of weighted {ENO} schemes.
\newblock {\em Journal of Computational Physics}, 126(1):202 -- 228, 1996.

\bibitem{Kanschat:2008}
Guido Kanschat.
\newblock Robust smoothers for high-order discontinuous {G}alerkin
  discretizations of advection–diffusion problems.
\newblock {\em Journal of Computational and Applied Mathematics}, 218(1):53 --
  60, 2008.

\bibitem{klofkorn2015efficient}
Robert Kl{\"o}fkorn.
\newblock Efficient matrix-free implementation of discontinuous {G}alerkin
  methods for compressible flow problems.
\newblock In {\em Proceedings of the Conference Algoritmy}, pages 11--21, 2015.

\bibitem{Krivodonova:2013}
Lilia Krivodonova and Ruibin Qin.
\newblock An analysis of the spectrum of the discontinuous {G}alerkin method.
\newblock {\em Applied Numerical Mathematics}, 64:1--18, 2013.

\bibitem{NME:NME5137}
M.~Kronbichler, S.~Schoeder, C.~M\"uller, and W.~A. Wall.
\newblock Comparison of implicit and explicit hybridizable discontinuous
  galerkin methods for the acoustic wave equation.
\newblock {\em International Journal for Numerical Methods in Engineering},
  106(9):712--739, 2016.
\newblock nme.5137.

\bibitem{VanLoan:2000}
Charles F.~Van Loan.
\newblock The ubiquitous {K}ronecker product.
\newblock {\em Journal of Computational and Applied Mathematics}, 123(1–2):85
  -- 100, 2000.
\newblock Numerical Analysis 2000. Vol. III: Linear Algebra.

\bibitem{Lynch:1964direct}
Robert~E. Lynch, John~R. Rice, and Donald~H. Thomas.
\newblock Direct solution of partial difference equations by tensor product
  methods.
\newblock {\em Numerische Mathematik}, 6(1):185--199, 1964.

\bibitem{Mardal:2007}
K.-A. Mardal, T.~K. Nilssen, and G.~A. Staff.
\newblock Order-optimal preconditioners for implicit {R}unge-{K}utta schemes
  applied to parabolic {PDEs}.
\newblock {\em {SIAM} Journal on Scientific Computing}, 29(1):361--375, 2007.

\bibitem{Orszag:1980}
Steven~A. Orszag.
\newblock Spectral methods for problems in complex geometries.
\newblock {\em Journal of Computational Physics}, 37(1):70 -- 92, 1980.

\bibitem{Peraire:2011}
Jaime Peraire and Per-Olof Persson.
\newblock High-order discontinuous {G}alerkin methods for {CFD}.
\newblock In Z.~J. Wang, editor, {\em Adaptive High-Order Methods in Fluid
  Dynamics}, chapter~5, pages 119--152. World Scientific, 2011.

\bibitem{Persson:2006efficient}
Per-Olof Persson and Jaime Peraire.
\newblock An efficient low memory implicit {DG} algorithm for time dependent
  problems.
\newblock In {\em 44th AIAA Aerospace Sciences Meeting and Exhibit}, page 113,
  2006.

\bibitem{Persson:2008}
Per-Olof Persson and Jaime Peraire.
\newblock {N}ewton-{GMRES} preconditioning for discontinuous {G}alerkin
  discretizations of the {N}avier-{S}tokes equations.
\newblock {\em SIAM Journal on Scientific Computing}, 30(6):2709--2733, 2008.

\bibitem{Reed:1973}
W.~H. Reed and T.~R. Hill.
\newblock Triangular mesh methods for the neutron transport equation.
\newblock {\em Los Alamos Report LA-UR-73-479}, 1973.

\bibitem{Saad:2003iterative}
Yousef Saad.
\newblock {\em Iterative methods for sparse linear systems}.
\newblock SIAM, 2003.

\bibitem{Shen2011}
Jie Shen, Tao Tang, and Li-Lian Wang.
\newblock {\em Separable Multi-Dimensional Domains}, pages 299--366.
\newblock Springer Berlin Heidelberg, Berlin, Heidelberg, 2011.

\bibitem{Shu2005}
Chi-Wang Shu, Wai-Sun Don, David Gottlieb, Oleg Schilling, and Leland Jameson.
\newblock Numerical convergence study of nearly incompressible, inviscid
  {T}aylor-{G}reen vortex flow.
\newblock {\em Journal of Scientific Computing}, 24(1):1--27, 2005.

\bibitem{Springel:2010}
Volker Springel.
\newblock E pur si muove: {G}alilean-invariant cosmological hydrodynamical
  simulations on a moving mesh.
\newblock {\em Monthly Notices of the Royal Astronomical Society}, 401(2):791,
  2010.

\bibitem{VanLoan:1993}
Charles~F. Van~Loan and Nikos Pitsianis.
\newblock Approximation with {K}ronecker products.
\newblock In {\em Linear Algebra for Large Scale and Real-Time Applications},
  pages 293--314. Springer, 1993.

\bibitem{van2011comparison}
Wim~M. Van~Rees, Anthony Leonard, D.~I. Pullin, and Petros Koumoutsakos.
\newblock A comparison of vortex and pseudo-spectral methods for the simulation
  of periodic vortical flows at high {R}eynolds numbers.
\newblock {\em Journal of Computational Physics}, 230(8):2794--2805, 2011.

\bibitem{Vos:2010}
Peter E.~J. Vos, Spencer~J. Sherwin, and Robert~M. Kirby.
\newblock From $h$ to $p$ efficiently: Implementing finite and spectral/$hp$
  element methods to achieve optimal performance for low-and high-order
  discretisations.
\newblock {\em Journal of Computational Physics}, 229(13):5161--5181, 2010.

\bibitem{Wang2013}
Z.J. Wang, Krzysztof Fidkowski, R{\'e}mi Abgrall, Francesco Bassi, Doru
  Caraeni, Andrew Cary, Herman Deconinck, Ralf Hartmann, Koen Hillewaert, H.T.
  Huynh, Norbert Kroll, Georg May, Per-Olof Persson, Bram van Leer, and Miguel
  Visbal.
\newblock High-order {CFD} methods: current status and perspective.
\newblock {\em International Journal for Numerical Methods in Fluids},
  72(8):811--845, 2013.

\bibitem{Warburton2008}
T.~Warburton and T.~Hagstrom.
\newblock Taming the {CFL} number for discontinuous galerkin methods on
  structured meshes.
\newblock {\em {SIAM} Journal on Numerical Analysis}, 46(6):3151--3180, jan
  2008.

\bibitem{yano2012optimization}
Masayuki Yano and David~L. Darmofal.
\newblock An optimization-based framework for anisotropic simplex mesh
  adaptation.
\newblock {\em Journal of Computational Physics}, 231(22):7626--7649, 2012.

\bibitem{Zahr2013}
Matthew~J. Zahr and Per-Olof Persson.
\newblock Performance tuning of {N}ewton-{GMRES} methods for discontinuous
  {G}alerkin discretizations of the {N}avier-{S}tokes equations.
\newblock In {\em 21st {AIAA} Computational Fluid Dynamics Conference}.
  American Institute of Aeronautics and Astronautics, June 2013.

\end{thebibliography}

\end{document}